%% file: subtorac.tex
\documentclass[12pt]{article}                                   

\usepackage{theorem}
\usepackage{graphics}
\usepackage{amsfonts}
\usepackage{amssymb}

\theoremstyle{change}{\theorembodyfont{\slshape}
\newtheorem{theorem}{Theorem.}[section]
\newtheorem{proposition}[theorem]{Proposition.}
\newtheorem{lemma}[theorem]{Lemma.}
\newtheorem{corollary}[theorem]{Corollary.}

 }

{
\theorembodyfont{\rmfamily}
\newtheorem{example}[theorem]{Example.}
\newtheorem{examples}[theorem]{Examples.}
\newtheorem{remark}[theorem]{Remark.}
\newtheorem{definition}[theorem]{Definition.}
}

\def\proof{\noindent {\bf Proof.}\enspace}
\def\endproof{ \quad $\kasten$}
\def\bemende{\quad \diamondsuit}

\def\S{{\cal S}}

\def\tpq{\mathop{\lower 4pt \hbox{$ {\buildrel{/} \over
        {\scriptstyle{\rm tpq}}} $}}\nolimits}

\def\mal{\! \cdot \!}
\def\h#1{\widehat{#1}}
\def\t#1{\widetilde{#1}}
\def\b#1{\overline{#1}}

\def\CC{{\mathbb C}}
\def\ZZ{{\mathbb Z}}
\def\RR{{\mathbb R}}

\def\PP{{\mathbb P}}
\def\mal{\mathbin{\! \cdot \!}}
\def\Orb{\mathop{\rm Orb}\nolimits}
\def\conv{\mathop{\hbox{\rm conv}}}
\def\cone{\mathop{\hbox{\rm cone}}}

\def\id{\mathop{\rm id}\nolimits}

\def\Hom{\mathop{\rm Hom}\nolimits}

\def\kasten{\mathord{\vbox{\hrule
                     \hbox{\vrule
                     \hskip5pt
                     \vrule height5pt
                     \vrule}
                     \hrule}}}
\def\text#1{\hbox{\rm #1}}

\def\bigtopmapleft#1{\smash{\mathop{\hbox to
      35pt{\leftarrowfill}}\limits^{#1}}}
\def\bigtopmapright#1{\smash{\mathop{\hbox to
      35pt{\rightarrowfill}}\limits^{#1}}}
\def\rmapdown#1{\Big\downarrow\rlap{$\vcenter
      {\hbox{$\scriptstyle#1$}}$}}
\def\lmapdown#1{\llap{$\vcenter{\hbox
      {$\scriptstyle#1$}}$}\Big\downarrow}

\begin{document}

\thispagestyle{empty}

\begin{center}
  
  {\LARGE\bf Toric Prevarieties
    
    \bigskip
    
    and Subtorus Actions}
  
  \bigskip
  
  \bigskip

Annette A'Campo--Neuen\footnote{email: {\tt Annette.ACampo@uni-konstanz.de}} %
and %
J\"urgen Hausen\footnote{email: {\tt Juergen.Hausen@uni-konstanz.de}}

\bigskip

Fachbereich Mathematik und Statistik

der Universit\"at Konstanz

\bigskip

\end{center}

%\bigskip

{\small
\noindent{\bf Abstract.}\enspace
Dropping separatedness in the definition of a toric variety, one
obtains the more general notion of a toric prevariety. Toric
prevarieties occur as ambient spaces in algebraic geometry and
moreover they appear naturally as intermediate steps in quotient
constructions. We first provide a complete description of the category
of toric prevarieties in terms of convex--geometrical data, so--called
systems of fans. In a second part, we consider actions of
subtori $H$ of the big torus of a toric prevariety $X$ and
investigate quotients for such actions. Using our language of systems
of fans,  we characterize existence of good prequotients for the
action of $H$ on $X$. Moreover, we show by means of an algorithmic
construction that there always exists a toric prequotient for the
action of $H$ on $X$, that means an $H$--invariant toric morphism $p$
from $X$ to a toric prevariety $Y$ such that every $H$--invariant
toric morphism from $X$ to a toric prevariety factors through $p$.
Finally, generalizing a result of D. Cox, we prove that every toric
prevariety $X$ occurs as the image of a categorical prequotient of 
an open toric subvariety of some $\CC^{s}$.
}

\section*{Introduction}

A {\it toric prevariety} is a normal prevariety, i.e.~possibly
non--separated, together with an effective regular action of an algebraic
torus that has a dense orbit. This notion occurs for the first
time in an article by J.~W\l{}odarczyk in 1991 (see \cite{Wlodar})
where he shows that  toric prevarieties are in fact universal ambient
spaces in algebraic geometry. More precisely, he proves: {\it Every normal
  variety over an algebraically closed field admits a closed embedding
  into a toric prevariety.}

In classical algebraic geometry, all the varieties considered
were quasi--projective; so
the ambient spaces were finite dimensional vector spaces or projective
spaces. But when in the 1940s  the abstract concept
of an algebraic variety via glueing affine charts was introduced, 
the resulting class of objects turned out to be much larger
than those fitting into the frame of the classical ambient spaces. It
therefore came as a surprise that the class of toric prevarieties is
indeed so large that they contain any given normal complex variety as
a closed subvariety. 

However, non--separated toric prevarieties have so far hardly
been studied. In the first part of this article our aim is to provide a
complete description of complex toric prevarieties and their morphisms in
terms of convex geometry. Generalizing the notion of a fan, we
introduce the concept of a {\it system of fans} in a lattice and obtain an
equivalence of the category of systems of fans and the category of
toric prevarieties (see Theorem \ref{equivcat}).

Our description makes the
category of toric prevarieties accessible for
explicit calculations. A first application is an
algorithmic construction of a {\it toric separation}: For every toric
prevariety $X$, we obtain a toric morphism from $X$ to a toric variety
$Y$ that is universal with respect to toric morphisms from $X$ to
toric varieties (see Theorem \ref{toricsep}). This result can serve
as a general tool for passing from the non--separated setting to
varieties. Further applications of the convex--geometrical language
are given in subsequent articles. For
example, it is used in \cite{ha} to prove a refined version of 
W\l odarczyk's embedding theorem.
 
The second part of this article is devoted to quotient constructions.
Frequently it is useful to decompose a quotient construction into a
non--separated first step followed by a separation process. For
example, A.~Bia\l ynicki--Birula uses non--separated quotient spaces
as a tool to prove in \cite{BB} that a normal variety contains only 
finitely many maximal subsets admitting a good quotient with respect
to a given reductive group action. 

In the toric setting it is natural to consider subtorus actions.
Here one faces the following problem:
Given a toric variety $X$ and a subtorus $H$ of the big torus of $X$;
when does there exist a suitable quotient for the action of $H$ on $X$?
This quotient problem has been studied by various authors: In
\cite{KaStZe} GIT--quotients for subtorus actions on projective 
toric varieties were investigated. More generally, J.~\'Swi\c{e}cicka
(\cite{Sw}) and H.~Hamm (\cite{Hm}) asked for existence of arbitrary
good quotients. In section~6 we treat their questions in the framework of 
non--separated toric prevarieties. In
analogy to the corresponding concept for the separated case, we 
define the notion of a {\it good prequotient}. We characterize in
terms of systems of fans when a good prequotient for the action of a
subtorus $H$ on a toric prevariety
$X$ exists (see Theorem \ref{characterization}). 

The characterizations given in the
above--mentioned results show that good quotients and prequotients
 only exist under quite special circumstances. So it is natural to ask
for more general notions. In \cite{acha}, the notion of a {\it toric
quotient} for a subtorus action on a toric variety was introduced and
it was proved that such a quotient  always exists. The analogous
notion in the context of toric prevarieties is the {\it toric
  prequotient}, i.e., an $H$--invariant toric morphism $p$ from $X$ to
a toric prevariety $Y$ 
such that every $H$--invariant toric morphism from $X$ factors
uniquely through $p$. In Section 7 we prove by means of an explicit
algorithm (see Theorem \ref{prequotcalc}) that toric prequotients
always exist. Determining first the toric prequotient of a subtorus
action on a toric variety and then performing toric separation splits
the calculation of the toric quotient into two steps (see Remark
\ref{intermediatestep}). This decomposition gives for example insight
into obstructions to the existence of categorical quotients (see
\cite{acha3}).

An application of our theory of quotients is given in the last
section: we represent an arbitrary toric prevariety as a quotient
space of an open toric subvariety of some $\CC^{n}$ by the action of a
subtorus of $(\CC^*)^n$ (see Corollary \ref{precox1}). This result
generalizes a similar statement on toric varieties due to D.~Cox (see
\cite{cox}). In contrast to Cox's construction, our quotient map is not
necessarily a good prequotient. However, it is  universal in the
category of algebraic prevarieties (see Proposition
\ref{precox1cat}). A toric prevariety $X$ with big torus $T$
occurs as the image of a good prequotient of an open subvariety of
some $\CC^{s}$ if and only if the intersection of any two maximal
affine open $T$--stable subspaces of $X$ is  affine (see Theorem
\ref{precox2suff}). 

\bigskip

\noindent %
We would like to thank A.~Bia\l ynicki--Birula, J.~Jurkiewicz,
J.~\'{S}wi\c{e}cicka and J.~W\l{}odarczyk for many suggestions and
helpful discussions. A particular thank goes to P. Heinzner for
pointing out a gap in an earlier version.

\section{Toric Prevarieties}

Let $X$ be a complex algebraic prevariety, i.e., a connected complex
ringed space that is obtained by glueing finitely many complex affine
varieties along open subspaces. Recall that $X$ is separated if and
only if it is Hausdorff with respect to the complex topology.

By definition a prevariety $X$ is normal if it is irreducible and all
its local rings are integrally closed domains, or equivalently, all
its affine charts are normal.  As in the case of varieties, a
normalization of a given prevariety $X$ is obtained by glueing
normalizations of affine charts of $X$.

\begin{definition}\label{torprevardef}
  A {\it toric prevariety} is a normal complex prevariety $X$ together
  with an effective regular action of an algebraic torus $T$ having an
  open orbit.
\end{definition}

For a toric prevariety $X$ and $T$ as in \ref{torprevardef}, we refer
to $T$ as the {\it acting torus} of $X$. Moreover, we fix a {\it base
  point} $x_{0}$ in the open orbit of $X$. Note that a toric
prevariety is a toric variety in the usual sense (see e.g. \cite{Fu})
if and only if it is separated.

\begin{example}\label{cdoubled0}
  The complex line, endowed with the $\CC^{*}$--action $t \mal z :=
  tz$, is a toric variety. Glueing two disjoint copies of $\CC$ along
  the open orbit $\CC^{*} \subset \CC$ yields a non-separated toric
  prevariety $X$. As a base point we choose $x_{0} := 1 \in \CC^{*}
  \subset X$. \quad $\diamondsuit$
\end{example}

\begin{proposition}\label{affcov}
  Every toric prevariety $X$ admits a finite covering of open affine
  subspaces that are stable by the acting torus $T$ of $X$.
\end{proposition}

\proof According to \cite{BB}, Theorem 1, there are only finitely many
maximal separated open subspaces $U_{i}$, $i \in I$, of $X$. Since $X$
is a noetherian topological space, it is covered by the
$U_{i}$.

By Sumihiro's Theorem (see \cite{Su}) we only have to show that each
$U_{i}$ is $T$--stable. This is done as follows: Every $t \in T$
permutes the sets $U_{i}$. Hence the elements of $T$ permute also the
complements $A_{i} := X \setminus U_{i}$. Consequently, for a given $i
\in I$ we have
$$ T = \bigcup_{j \in I} {\rm Tran}_{T}(A_{i},A_{j}),$$
where ${\rm Tran}_{T}(A_{i},A_{j})$ denotes the closed set $\{t \in T;
\; t \mal A_{i} \subset A_{j}\}$. Since $T$ is irreducible, there is a
$j_{0} \in I$ such that $T = {\rm Tran}_{T}(A_{i},A_{j_{0}})$. Note
that $A_{i} = e_{T} \mal A_{i} \subset A_{j_0}$. Thus maximality of
$U_{j_{0}}$ implies $i = j_{0}$ which yields $T \mal U_{i} =
U_{i}$. \endproof

\bigskip

Note that the arguments of the above proof yield that any
$G$--prevariety, where $G$ is a connected algebraic group, can be
covered by $G$--stable separated open subspaces. For disconnected $G$
this statement is false (see Example 1.6).

Now assume that $X$ is a toric prevariety with acting torus $T$. Using
the theory of toric varieties we can conclude from Proposition \ref{affcov}
that the set $\Orb(X)$ of all $T$--orbits of $X$ is
finite. %Furthermore, we have

\begin{remark}\label{affdecomp}
  For every $T$--orbit $B$ of $X$ there exists a unique $T$-stable
  open affine subspace $X_{B}$ of $X$ such that $B$ is a closed subset
  of $X_{B}$. Moreover, we have
  $$
  X_{B} = \bigcup_{B' \in \Orb(X); \; B \subset \b{B'}} B' . 
\qquad\bemende
%\leqno{(*)}
$$
\end{remark}

%
%\proof Proposition \ref{affcov} and the theory of  toric varieties
%yield that every $T$--orbit $B$ has an affine $T$-stable neighbourhood $X_{B}$
%such that $B$ is closed in $X_{B}$. In order to prove the assertion we
%only have to verify the equation $(*)$.
%
%The inclusion ``$\subset$'' is a basic observation in the theory of
%toric varieties (see e.g. \cite{Fu}, Section 3.1). In order to verify
%``$\supset$'', consider a $T$--orbit $B' \subset X \setminus
%X_{B}$. Since $X_{B}$ is open, we have $\b{B'} \subset X \setminus
%X_{B}$ and, in particular, $B \not\subset \b{B'}$.  \endproof
%

The morphisms in the category of toric prevarieties are defined
similarly as in the separated case: Let $f \colon X \to X'$ be a
regular map of toric prevarieties $X$ and $X'$ with base points
$x_{0}$ and $x_{0}'$ respectively.

\begin{definition}
  The map $f$ is called a {\it toric morphism} if $f(x_{0}) = x'_{0}$
  and there is a homomorphism $\varphi \colon T \to T'$ of the acting
  tori of $X$ and $X'$ such that $f(t \mal x) = \varphi(t) \mal f(x)$
  holds for all $(t,x) \in T \times X$.
\end{definition}

\begin{example}
  For the toric prevariety $X$ of Example \ref{cdoubled0} let $0_{1}$
  and $0_{2}$ denote the two fixed points of $T$. Then
  $f\vert_{\CC^{*}} := \id_{\CC^{*}}$, $f(0_{1}) := 0_{2}$ and
  $f(0_{2}) := 0_{1}$ defines a toric automorphism of $X$ of order $2$.
\end{example}

\begin{lemma}\label{mappingaffinecharts}
  Let $f\colon X \to X'$ be a toric morphism of toric prevarieties and
  let $B \subset X$, $B' \subset X'$ be orbits of $T$ and $T'$
  respectively. Then we have $f(X_{B}) \subset X'_{B'}$ if and only if
  $f(B) \subset X'_{B'}$.
\end{lemma}

\proof Suppose that $f(B) \subset X'_{B'}$ holds. Then Remark
\ref{affdecomp} yields $B'\subset \b{T'\mal f(B)}$. Now consider a
$T$--orbit $B_{1}$ with $B \subset \b{B_{1}}$. Then $f(B)\subset
\b{f(B_{1})}$ and hence $B'$ is contained in the closure of the
$T'$-orbit $T'\mal f(B_{1})$. That means that $f(B_{1}) \subset
X'_{B'}$. Thus Remark \ref{affdecomp} implies $f(B) \subset
X'_{B'}$. \endproof

\section{Systems of Fans}

In this section we introduce the notion of a system of fans in a
lattice and associate to every system of fans a toric prevariety. Our
construction is a generalization of the basic construction in the
theory of toric varieties. First we have to fix some notation:

By a lattice we mean a free $\ZZ$--module of finite rank. For a given
lattice $N$ let $N_{\RR} := \RR \otimes_{\ZZ} N$ denote its associated
real vector space. Moreover, for a homomorphism $F \colon N \to N'$ of
lattices, denote by $F_{\RR}$ its extension to the real vector spaces
associated to $N$ and $N'$.

In the sequel let $N$ be a lattice. When we speak of a cone in $N$ we
always think of a (not necessarily strictly) convex rational
polyhedral cone in $N_{\RR}$. For a cone $\sigma$ in $N$ we denote by
$\sigma^{\circ}$ the relative interior of $\sigma$ and if $\tau$ is a
face of $\sigma$, then we write $\tau \prec \sigma$.

As usual, we call a finite set $\Delta$ of strictly convex cones in
$N$ a {\it fan\/} in $N$ if any two cones of $\Delta$ intersect in a
common face and if $\sigma \in \Delta$ implies that also every face of
$\sigma$ lies in $\Delta$. If $\Delta'$ is a subfan of a fan $\Delta$
we will write $\Delta' \prec \Delta$. A fan $\Delta$ is called {\it
  irreducible\/} if it consists of all the faces of a cone $\sigma$.

\begin{definition}\label{systemoffans}
  Let $I$ be a finite index set. A collection
  $\S=(\Delta_{ij})_{i,j\in I}$ of fans in $N$ is called a {\it system
    of fans\/} if the following properties are satisfied for all
  $i,j,k \in I$:
\begin{enumerate}
\item $\Delta_{ij} = \Delta_{ji}$,
\item $\Delta_{ij} \cap \Delta_{jk} \prec \Delta_{ik}$.
\end{enumerate}
Note that in particular, $\Delta_{ij} \prec \Delta_{ii} \cap
\Delta_{jj}$
for all $i,j\in I$.
\end{definition}

\begin{examples}\label{fan2prefan}
\begin{enumerate}
\item Every fan $\Delta$ in $N$ can be considered as a system of fans
  $\S = (\Delta)$ with just one element.
\item Let $\sigma_{1}, \ldots, \sigma_{r}$ denote the maximal cones of
  a fan $\Delta$ and set $I := \{1, \ldots, r\}$.  Let $\Delta_{ii}$
  denote the fan of faces of $\sigma_i$ and define
  $\Delta_{ij}:=\Delta_{ii}\cap\Delta_{jj}$. Then $\S :=
  (\Delta_{ij})_{i,j\in I}$ is a system of irreducible fans.
\item If $\Delta_{1}, \ldots, \Delta_{r}$ are fans in a lattice $N$,
  then $\Delta_{ii} := \Delta_{i}$ and $\Delta_{ij} := \{ \{ 0 \} \}$
  defines a system of fans in $N$.
\item For a given collection $\sigma_{1}, \ldots, \sigma_{r}$ of
  strictly convex cones in a lattice $N$, set $I := \{1, \ldots, r\}$.
  Let $\Delta_{ii}$ denote the fan of faces of $\sigma_i$ and define
  $\Delta_{ij}$ to be the fan of all common proper faces of
  $\sigma_{i}$ and $\sigma_{j}$. Then $\S := (\Delta_{ij})_{i,j\in I}$
  is a system of fans. \quad $\diamondsuit$
\end{enumerate}
\end{examples}

In the sequel let $\S=(\Delta_{ij})_{i,j\in I}$ be a given system of
fans in $N$. There is an algebraic torus having $N$ as its lattice of
one parameter subgroups, namely $T:=\Hom(N^\vee,\CC^*)$, where $N^\vee
:= \Hom(N,\ZZ)$ denotes the dual module of $N$. We associate to $\S$ a
toric prevariety $X_{\S}$ with acting torus $T$ as follows:
 
For each index $i\in I$ let $X_{i} := X_{\Delta_{ii}}$ denote the
toric variety corresponding to the fan $\Delta_{ii}$ (see e.g.
\cite{Fu}).  For any two indices $i \ne j$ let $X_{ij}$ and $X_{ji}$
be the open toric subvarieties of $X_{i}$ and $X_{j}$ corresponding to
the subfan $\Delta_{ij}$.

The lattice homomorphism $\id_{N}$ defines toric isomorphisms $f_{ji}
\colon X_{ij} \to X_{ji}$. Note that Property \ref{systemoffans} ii)
yields $f_{kj} \circ f_{ji} = f_{ki}$ on the intersections $X_{ij}
\cap X_{ik}$. Define $X_{\S}$ to be the $T$--equivariant glueing of
the toric varieties $X_{i}$ by the glueing maps $f_{ij}$. By
construction, $X_{\S}$ is a toric prevariety.

\begin{example}
  Let $I := \{1,2\}$ and let $\Delta_{11} = \Delta_{22}:= \{ \{0\},
  \sigma \}$ be the fan of faces of $\sigma := \RR_{\ge 0}$.  Setting
  $\Delta_{12} := \Delta_{21} := \{ \{0\}\}$ we obtain a system of
  fans $\S$ in $\ZZ$.  The associated toric prevariety $X_{\S}$ is the
  complex line with zero doubled (see Example \ref{cdoubled0}). \quad
  $\diamondsuit$
\end{example}

Apparently, different systems $\S$ and $\S'$ of fans can lead to the
same toric prevariety $X$, since there may be various possibilites for
choosing separated toric charts covering the prevariety $X$.

\begin{example}
  Let $N := \ZZ$ and $\sigma := \RR_{\ge 0}$. The fan $\Delta:=
  \{\sigma,-\sigma, \{0\} \}$ gives rise to the toric variety
  $X_{\Delta} = \PP_1$. Let $I := \{1,2\}$ and set
  $$\Delta_{11} := \Delta_{22} := \Delta, \qquad \Delta_{12}
  :=\Delta_{21} := \{ \{0\}\}.$$
  Then the resulting system of fans $\S
  := (\Delta_{ij})_{i,j\in I}$ defines the toric prevariety $X$ that
  is obtained from glueing two copies of $\PP_1$ along $\CC^*$. If we
  set $I' := \{1,2,3,4\}$,
  $$
  \Delta'_{11} := \Delta'_{22}:= \{ \{0\}, \sigma \}, \qquad
  \Delta'_{33} := \Delta'_{44} := \{ \{0\}, -\sigma \}, $$
  and
  $\Delta'_{ij} := \{ \{0\} \}$ for all $i \ne j$, then we arrive at a
  system $\S'$ of fans defining the same toric prevariety $X$ as
  above. But now the fans of the system are irreducible and correspond
  to affine charts of $X$.  \quad $\diamondsuit$
\end{example}

A given toric prevariety has two distinguished systems of charts,
namely the covering by maximal $T$--stable separated charts and the
covering by maximal $T$--stable affine charts. The latter one
corresponds to systems $\S = (\Delta_{ij})_{i,j\in I}$ of fans where
every $\Delta_{ii}$ is irreducible. Such a system will be called {\it
  affine\/}.

For the description of the orbit structure of $X_{\S}$ the following
observation will turn out to be useful: The system of fans $\S$
naturally induces an equivalence relation on the set $\mathfrak{F}(\S)
:= \{(\sigma,i); \; i \in I,\sigma \in \Delta_{ii}\}$ of labelled
cones, namely
$$(\sigma,i) \sim (\sigma,j) \quad \Longleftrightarrow \quad
\sigma\in\Delta_{ij}\,.$$
We call this equivalence relation the {\it
  glueing relation\/} of $\S$, and we denote the set of equivalence
classes by $\Omega:=\Omega(\S)$. The equivalence class of an element
$(\sigma, i) \in \mathfrak{F}(\S)$ is denoted by $[\sigma,i]$.

\begin{remark}\label{glueingrelation}
  The glueing relation satisfies the following conditions:
\begin{enumerate}
\item $(0, i) \sim (0, j)$ for all $i$, $j$,
\item $(\sigma,i) \sim (\tau,j)$ implies $\sigma = \tau$,
\item $(\tau,i) \sim (\tau,j)$ implies $(\tau',i) \sim (\tau',j)$ for
  every $\tau' \prec \tau$. $\bemende$
\end{enumerate}
\end{remark}

As a converse of Remark \ref{glueingrelation}, we can recover $\S$
from its glueing relation: Let $S$ denote a finite set of cones in
$N$, let $I$ be a finite index set and suppose that $\mathfrak{F}$ is
a subset of $S\times I$ where for every $i$ the set $\Delta_{ii}:=
\mathfrak{F} \cap (S\times \{i\})$ forms a fan.

\begin{remark}
  If $\sim$ is an equivalence relation on $\mathfrak{F}$ satisfying
  the conditions \ref{glueingrelation} i)--iv), then we obtain a
  system of fans by setting
  $$\Delta_{ij} := \{\tau\in\Delta_{ii}\cap\Delta_{jj} ; \; (\tau,i)
  \sim (\tau,j)\}\,.\quad\bemende
  $$
\end{remark}

Let us now return to the toric prevariety $X_{\S}$ obtained from
glueing the toric varieties $X_i=X_{\Delta_{ii}}$. Recall that the
$T$--orbits of $X_i$ are in $1$--$1$--correspondence with the cones in
$\Delta_{ii}$. For every $\sigma\in\Delta_{ii}$, there is even a
distinguished point $x_{(\sigma,i)}$ in the corresponding $T$--orbit
in $X_{i}$ (see e.g.~\cite{Fu}, p.28).

In the toric prevariety $X_{\S}$ a point $x_{(\sigma,i)}\in X_{i}$ is
identified with $x_{(\tau,j)} \in X_{j}$ if and only if
$x_{(\sigma,i)} \in X_{ij}$ and $x_{(\tau,j)} \in X_{ji}$ and $\sigma
= \tau$, or equivalently if $(\sigma,i) \sim (\tau,j)$.  So a
distinguished point $x_{(\sigma,i)} \in X_{i}$ defines a {\it
  distinguished point} in $X_{\S}$ which depends only on the
equivalence class $[\sigma,i]$ of $(\sigma,i)$ in $\Omega(\S)$ and is
denoted by $x_{[\sigma,i]}$.

\begin{remark}\label{distpoints}
  The assignment $[\sigma,i] \mapsto T \mal x_{[\sigma,i]}$ defines a
  bijection from $\Omega(\S)$ to the set of $T$-orbits of the toric
  prevariety $X_{\S}$.  $\bemende$
\end{remark}

The point $x_{0} := x_{[\{0\},i]}$ corresponding to the open $T$-orbit
will be considered as the {\it base point} of $X_{\S}$.  For a
distinguished point $x_{[\sigma,i]}$ of $X_{\S}$ we define
$X_{[\sigma,i]}$ to be the open affine $T$--stable subspace of
$X_{\S}$ that contains $T \mal x_{[\sigma,i]}$ as closed subset (see
Remark 1.4).
  
By Property \ref{glueingrelation} iv), the face relation induces a
partial ordering on the set ${\Omega}(S)$, namely $[\tau,j] \prec
[\sigma,i]$ if $\tau$ is a face of $\sigma$ and $[\tau,i] = [\tau,j]$.
This partial ordering reflects the behaviour of orbit closures in
$X_{\S}$:

\begin{lemma}\label{orbitclosures}
  A point $x_{[\sigma, i]}$ lies in the closure of the orbit $T \mal
  x_{[\tau,j]}$ if and only if $[\tau, j] \prec [\sigma, i]$. In
  particular, one has
  $$
  X_{[\sigma,i]} = \bigcup_{[\tau, j] \prec [\sigma, i]} T \mal
  x_{[\tau,j]} = \bigcup_{[\tau, j] \prec [\sigma, i]} X_{[\tau,j]}.
  $$
\end{lemma}

\proof Assume $[\tau, j] \prec [\sigma, i]$. By definition of the
partial ordering ``$\prec$'', this means $\tau \prec \sigma$ and
$[\tau,j] = [\tau,i]$. This implies $x_{(\sigma, i)} \in \b{T \mal
  x_{(\tau,i)}} \subset X_{i}$. Hence $x_{[\sigma,i]}$ lies in the
closure of $T \mal x_{[\tau, j]}$.

Now, let $x_{[\sigma,i]} \in \b{T \mal x_{[\tau, j]}}$. Consider the
$T$--stable separated neighbourhood $X_i$ of $x_{[\sigma,i]}$. Since
$X \backslash X_i$ is closed, we have $x_{[\tau,j]} \in X_i$, i.e.,
$[\tau,j] = [\tau,i]$. Now the theory of toric varieties tells us that
in $X_i$ we have $\tau \prec \sigma$.  \endproof

Together with the corresponding statement on affine toric varieties,
Lemma \ref{mappingaffinecharts} implies the following

\begin{remark}\label{fmapsdist}
  Let $f \colon X_{\S} \to X_{\S'}$ be a toric morphism. Then $f$ maps
  distinguished points to distinguished points.  $\bemende$
\end{remark}

\section{Toric Morphisms and Maps of Systems of Fans}

We first introduce the concept of a map of systems of fans and then show that
$\S \mapsto X_{\S}$ is an equivalence of categories. Let
$\S=(\Delta_{ij})_{i,j\in I}$ and $\S'=(\Delta_{ij}')_{i,j\in I'}$ denote
systems of fans in lattices $N$ and $N'$ respectively.

\begin{definition}\label{mapofsystemsoffans}
  A {\it map of systems of fans} from $\S$ to $\S'$ is a pair
  $(F,\mathfrak{f} )$, where $F \colon N \to N'$ is a lattice
  homomorphism and $\mathfrak{f} \colon \Omega(S) \to \Omega(S')$ is
  a map with the following properties:
\begin{enumerate}
\item If $[\tau,j] \prec [\sigma,i]$ then $\mathfrak{f}([\tau,j])
  \prec \mathfrak{f}([\sigma,i])$, i.e., $\mathfrak{f}$ is order
  preserving.
\item If $\mathfrak{f}([\sigma,i]) = [\sigma',i']$ then
  $F_{\RR}(\sigma^{\circ}) \subset (\sigma')^{\circ}$.
\end{enumerate}
\end{definition}

\begin{remark}\label{maps2complete}
  Assume that $\S'=(\Delta')$ is a single fan in $N'$ and $F \colon N\to N'$
  is a lattice homomorphism such that $F_{\RR}$ maps the cones of $\S$
  into cones of $\Delta'$. Then there is a unique map $\mathfrak{f} \colon
  \Omega(\S) \to {\Omega}(\S')$ such that 
  $(F,\mathfrak{f})$ is a map of
  the systems of fans $\S$ and $\S'$. $\bemende$
\end{remark}

On the other hand, if $\S'$ is arbitrary but $\S=(\Delta)$ is a single
fan and $F \colon N \to N'$ is a lattice homomorphism such that
$F_{\RR}$ maps cones of $\S$ into cones of $\S'$, then there need not
exist a map $(F,\mathfrak{f})$ of the systems of fans $\S$ and $\S'$,
even if the glueing relation of $\S'$ is maximal:

\begin{example}
  Let $\Delta$ be the fan in $\ZZ^{3}$ having $\sigma_{1} :=
  \cone(e_1,e_2,e_1+e_2+e_3)$ and $\sigma_{2} :=
  \cone(e_2,e_1+e_2+e_3,e_2-e_1+e_3)$ as its maximal cones.  Let $F
  \colon \ZZ^{3} \to \ZZ^{2}$ denote the projection given by
  $(x,y,z)\mapsto (x,y)$.
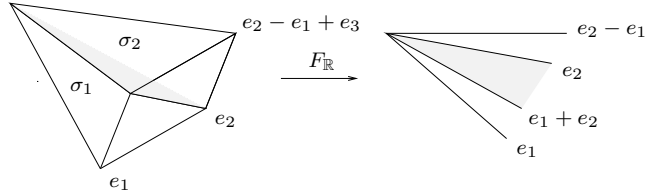
\begin{center}
  \input{nomap.pstex_t}
\end{center}
Let $I':=\{1,2\}$ and let $\Delta_{ii}$ be the fan of faces of the cones
$\tau_i:=F_{\RR}(\sigma_i)$ in $\ZZ^{2}$. Let $\S'$ denote the system
of fans obtained obtained from the $\Delta_{ii}$ by adding $\Delta_{12}:=0$. 
Then there is no map of systems of fans $(F,\mathfrak{f})$ from
$\S$ to $\S'$. \quad $\diamondsuit$
\end{example}

To obtain a functor from the category of systems of fans to the
category of toric prevarieties we now define the assignement on the
level of morphisms. Let $X := X_{\S}$ and $X' := X_{\S'}$ be the
respective toric prevarieties arising from $\S$ and $\S'$ and let
$(F,\mathfrak{f})$ be a map of the systems of fans $\S$ and $\S'$.

We construct a toric morphism $f \colon X \to X'$ as follows. 
For a given $i \in I$ and $\sigma\in\Delta_{ii}$, 
set $[\sigma',i'] := \mathfrak{f}([\sigma,i])$.
Then $F_{\RR}(\sigma) \subset \sigma'$, so $F$ defines a toric
morphism $f_{[\sigma,i]}$ from the affine toric variety 
$X_{[\sigma,i]}$ to $X'_{[\sigma',i']}$. By condition
\ref{mapofsystemsoffans} i) we obtain
$$ f_{[\sigma,i]} \vert_{X_{[\sigma,i]}\cap X_{[\tau,j]}} 
=f_{[\tau,j]} \vert_{X_{[\sigma,i]}\cap X_{[\tau,j]}}$$
for every $[\tau,j] \in \Omega(\S)$. Consequently the regular maps
$f_{[\sigma,i]}$ glue together to a toric morphism $f \colon X \to X'$.

The geometric meaning of the map $\mathfrak{f} \colon \Omega(\S) \to
\Omega(\S')$ is to prescribe the values of the distinguished points
for the toric morphism $f$ associated to a map $(F, \mathfrak{f})$ of
the systems $\S$ and $\S'$ of fans:

\begin{lemma}\label{nufgeom}
  For every $[\sigma,i] \in \Omega(\S)$ we have $f(x_{[\sigma, i]}) =
  x'_{\mathfrak{f}([\sigma, i])}$.
\end{lemma}
  
\proof Let $[\sigma,i] \in \Omega(\S)$ and let $[\sigma',i'] :=
\mathfrak{f}([\sigma,i])$. Consider the toric morphism $f_{i} := f_{[\sigma,i]}
\colon X_{[\sigma,i]} \to X'_{[\sigma',i']}$. Since
$F_{\RR}(\sigma^{\circ})$ is contained in $(\sigma')^{\circ}$, we have
$$ f_{i}(x_{(\sigma,i)}) = f_{i}\left(\lim_{t \to 0}
\lambda_{v}(t) \mal x_{(\{0\},i)}\right) = \lim_{t \to 0}
f_{i}(\lambda_{v}(t) \mal x_{(\{0\},i)}) = \lim_{t \to 0}
\lambda_{F(v)}(t) \mal x_{(\{0\},i)} = x_{(\sigma',i')}$$
where $v$ is any lattice point in $\sigma^{\circ}$ and $\lambda_{v}
\colon \CC^{*} \to T$ denotes the one--parameter--subgroup of the
acting torus $T$ of $X_{[\sigma,i]}$ defined by $v$. This yields the
claim. \endproof

\bigskip

Denote by $\varphi \colon T \to T'$ the homomorphism of acting tori
associated to the toric morphism $f$. Then we obtain the following
description of the fibers of $f$.

\begin{proposition}\label{fibreformula}
  {\it Fibre Formula.}\enspace For every distinguished point
  $x_{[\sigma',i']} \in X_{\S'}$ we have
$$ f^{-1}(x_{[\sigma',i']}) \; = \; %
\bigcup\! \! \lower10pt\hbox{$\scriptstyle%
                         [\sigma,i]\in\mathfrak{f}^{-1}([\sigma',i'])$ }%
                       \varphi^{-1}(T'_{x_{[\sigma',i']}}) \mal
                       x_{[\sigma,i]} .
                       $$
\end{proposition}

\proof The inclusion ``$\supset$'' follows from Lemma
\ref{nufgeom}. In order to check ``$\subset$'', let $x \in
f^{-1}(x_{[\sigma',i']})$. Then $x = t \mal x_{[\sigma,i]}$ for some
$t \in T$ and $[\sigma,i] \in \Omega(\S)$ and hence
$$f(x) = x_{[\sigma',i']} = \varphi(t) \mal f(x_{[\sigma,i]}).$$
By Lemma \ref{nufgeom}, we know obtain that $f(x_{[\sigma,i]})$ is a
distinguished point. This implies 
$$f(x_{[\sigma,i]}) = x_{[\sigma',i']} \quad \hbox{and} \quad
\varphi(t) \mal x_{[\sigma',i']} = x_{[\sigma',i']}. \quad \kasten $$

\bigskip

Now we come to the main result of this section, namely to generalize
the correspondence between fans and toric varieties to a
correspondence between systems of fans and toric prevarieties. By
construction, $\mathfrak{TP} \colon \S \mapsto X_{\S}$,
$(F,\mathfrak{f}) \mapsto f$ is a covariant functor from the category
of systems of fans to the category of toric prevarieties.

\begin{theorem}\label{equivcat}
  $\mathfrak{TP}$ and the restriction of $\mathfrak{TP}$ to the
  (full) subcategory of affine systems of fans are equivalences
  of categories.
\end{theorem}

\proof By equivariance, a toric morphism is determined by its associated
homomorphism of the acting tori and its values on the distinguished
points. Hence Lemma \ref{nufgeom} yields that the functor $\mathfrak{TP}$ is
faithful.

Next we verify that $\mathfrak{TP}$ is fully faithful. Let $\S$
and $\S'$ be systems of fans in lattices $N$ and $N'$ respectively and let
$f \colon X_{\S} \to X_{\S'}$ be a toric morphism.  Then the
associated homomorphism $\varphi \colon T \to T'$ of the respective
tori defines a homomorphism $F \colon N \to N'$.

For $[\sigma,i] \in \Omega(\S)$ the associated distinguished point
$x_{[\sigma,i]}$ is mapped to a distinguished point $x'_{[\tau,j]}$
and $f(X_{[\sigma,i]}) \subset X'_{[\tau,j]}$ (see Lemma
\ref{mappingaffinecharts} and Remark \ref{fmapsdist}).  Set
$\mathfrak{f}([\sigma,i]) := [\tau,j]$. Now it follows from Remark
\ref{affdecomp} and Lemma \ref{orbitclosures} that $(F,\mathfrak{f})$ is a
map of systems of fans. By Lemma \ref{nufgeom}, $f$ is the toric morphism
associated to $(F,\mathfrak{f})$.

Finally we have to show that for every toric prevariety $X$ there
exists an affine system of fans $\S$ with $X \cong X_{\S}$. 
Let $X_1,\dots,X_r$ be the maximal $T$-stable affine charts of $X$
(see Proposition \ref{affcov}), and let $N$ denote the
lattice of one--parameter--subgroups of the acting torus $T$ of $X$.
Then every chart $X_i$ corresponds to a irreducible fan $\Delta_{ii}$
in $N$.

For every $i,j\in I:=\{1,\dots,r\}$ the  intersection $X_i\cap X_j$ 
is  an open toric subvariety of both $X_i$ and $X_j$ 
and hence corresponds to a fan $\Delta_{ij}$
which is a common subfan of $\Delta_{ii}$ and $\Delta_{jj}$.
It follows that the collection $\S:=(\Delta_{ij})_{1 \le i, j \le r}$
forms an affine system of fans and it is straightforward to check that
$X \cong X_{\S}$.  \endproof

\section{The Toric Separation}

Let $X$ be a toric prevariety. A {\it toric separation\/} of $X$ is a
toric morphism $p \colon X \to Y$ to a toric variety $Y$ that has the
following universal property: For every toric morphism $f$ from $X$ to
a toric variety $Z$ there exists a unique toric morphism $\t{f} \colon
Y \to Z$ such that $f = \t{f} \circ p$ holds. The main result of this
section is

\begin{theorem}\label{toricsep}
  Every toric prevariety has a toric separation.
\end{theorem}

We prove this statement by showing the corresponding result (Theorem
\ref{red2fan} below) in the category of systems of fans. First we
translate the notion of the toric separation into the language of
systems of fans: Let $N$ be a lattice and assume that $\S$ is a system
of fans in $N$.

We call a map $(P,\mathfrak{p})$ of systems of fans from $\S$ to a fan
$\t{\Delta}$ in a lattice $\t{N}$ a {\it reduction to a fan}, if for
each map $(F,\mathfrak{f})$ of systems of fans from $\S$ to a fan
$\Delta'$ in a lattice $N'$ there is a unique lattice homomorphism
$\t{F}\colon \t{N}\to N'$ defining a map of the fans $\t{\Delta}$ and
$\Delta'$ such that $F = \t{F} \circ P$.

\begin{theorem}\label{red2fan}
  Every system of fans admits a reduction to a fan.
\end{theorem}

\proof Let $\S$ be a system of fans in a lattice $N$. Then
$S:=\bigcup_i \Delta_{ii}$ is a system of cones in $N$ in the sense of
\cite{acha}, Section 2. Moreover, every map of systems of fans from
$\S$ to a fan $\Delta'$ is also a map of the systems $S$ and $\Delta'$
of cones.

Denote by $\t{\Delta}$ the quotient fan of $S$ by the trivial
sublattice $L= \{0\}$ of $N$ (see \cite{acha}, Definition 2.1 and
Theorem 2.3). Then $\t{\Delta}$ lives in some lattice $\t{N}$ and
there is a projection $P \colon N \to \t{N}$, mapping the cones from
$S$ into cones of $\t{\Delta}$. By Remark \ref{maps2complete}, $P$
defines a map of the systems of fans $\S$ and $\t{\Delta}$. It follows
from the universal property of the quotient fan that $P$ is the
reduction of $\S$ to a fan. \endproof

\bigskip

By a {\it separation} of a prevariety $X$ we mean a regular map $p$
from $X$ to a variety $Y$ that is universal with respect to arbitrary
regular maps from $X$ to varieties. It can be shown that every toric
prevariety of dimension less than three has such a separation (see
\cite{acha3}).

In dimension three we find the first examples of toric prevarieties
that need not have a separation. The remainder of this section is
devoted to giving such an example. Let $I=\{1,2\}$ and let $\S$ be the
affine system of fans in $\ZZ^{3}$ determined by the cones
$$\sigma_1 := \cone(e_{1},e_{2}), \qquad \sigma_2 :=
\cone(e_{1}+e_{2}, e_{3})$$
glued along $0$. Note that $X_{\S}$ is the
glueing of $\CC^{2} \times \CC^{*}$ and $\CC^{*} \times \CC^{2}$ along
$(\CC^{*})^{3}$ via the map $(t_{1}, t_{2}, t_{3}) \mapsto
(t_{1}t_{2}^{-1}, t_{2}, t_{3})$.

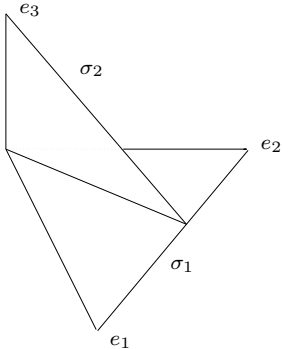
\begin{center}
  \input{nonsurj.pstex_t}

\end{center}

\begin{proposition}
  $X_{\S}$ admits no separation.
\end{proposition}

\proof Assume that there exists a separation $p \colon X_{\S} \to Y$.
With the universal property we obtain that $p$ is surjective, $Y$ is
normal and there is an induced (set theoretical) action of $T :=
(\CC^{*})^{3}$ on $Y$ such that $p$ is equivariant. We lead this to a
contradiction by showing that the toric separation of $X_\S$ does not
factor through $p$.

First we describe the toric separation of $X_\S$ explicitly. Let
$\t{\Delta}$ be the fan of faces of $\sigma :=
\cone(e_{1},e_{2},e_{3})$ in $\ZZ^{3}$. The reduction of $\S$ to a fan
is the map $Q := \id_{\ZZ^{3}}$ of the systems of fans $\S$ and
$\t{\Delta}$. Set $Z := X_{\t{\Delta}} = \CC^{3}$. The toric
separation of $X_{\S}$ is the toric morphism $q \colon X_{\S} \to Z$
associated to the map $Q$ of systems of fans. Note that $q(X_{\S})$ is
not open in $Z$, since we have
$$q(X_{\S}) = \CC^{3} \setminus (\{0\} \times \CC^{*} \times \{0\}
\cup \CC^{*} \times \{0\} \times \{0\}). $$
Now, by the universal property of $p$ there is a unique regular map $f
\colon Y \to Z$ such that $q = f \circ p$. Clearly $f$ is
$T$-equivariant. We claim moreover that $f$ is injective. To verify
this, we investigate the fibres of $f$.

Note first that, by equivariance of $f$, it suffices to consider the
fibres of distinguished points. Moreover, by surjectivity of $p$, we
have $f^{-1}(z) = p(q^{-1}(z))$ for every $z \in Z$.  Using the Fibre
Formula \ref{fibreformula}, we see that
$$
q^{-1}(z_{0}) = x_{0}, \qquad q^{-1}(z_{\varrho_{i}}) =
x_{\varrho_{i}},$$
where $\varrho_{i} := \RR_{\ge 0} e_{i}$. This implies that
$f^{-1}(z_{0})$ as well as the fibres $f^{-1}(z_{\varrho_{i}})$
consist of exactly one point. Again by the Fibre Formula one has
$$
q^{-1}(z_{\sigma_1}) = \{ x_{[\sigma_1,1]} \} \cup H \mal
x_{[\varrho,2]},$$
where $\varrho := \RR_{\ge 0} (e_{1}+e_{2})$ and $H$ is the subtorus
of $T$ that corresponds to the sublattice $\ZZ e_{1} \oplus \ZZ e_{2}$
of $\ZZ^{3}$. Note that for the one parameter subgroup $\lambda \colon
\CC^{*} \to T$ corresponding to the lattice vector $e_{1} + e_{2}$ we
obtain
$$
\lim_{t \to 0} \lambda(t) \mal x_0 = x_{[\sigma_1,1]} , \qquad
\lim_{t \to 0} \lambda(t) \mal x_0 = x_{[\varrho,2]}$$
in the affine
charts $X_{\sigma_1}$ and $X_{\sigma_2}$ respectively. Thus the points
$x_{[\sigma_1,1]}$ and $x_{[\varrho,2]}$ cannot be separated by
complex open neighbourhoods.

Since $p$ is continuous with respect to the complex topology and $Y$
is Hausdorff it follows $p(x_{[\sigma_1,1]}) = p(x_{[\varrho,2]})$.
Since $H$ fixes $x_{[\sigma_1,1]}$ and $p$ is equivariant, $H$ fixes
also $p(x_{[\varrho,2]})$. Consequently we obtain that
$f^{-1}(z_{\sigma_1})$ consists of a single point.

Finally we have to consider $z_{\sigma}$. We have $q^{-1}(z_{\sigma})
= T \mal x_{[\sigma_2, 2]}$. In order to see that $f^{-1}(z_{\sigma})$
is a single point, it suffices to check that $p(x_{[\sigma_2, 2]})$ is
fixed by $T$. Since $q(x_{[\sigma_2,2]}) \ne q(x_{[\varrho,2]})$, we
obtain
\begin{eqnarray*}
T \mal p(x_{[\sigma_2, 2]}) & \subset & \b{p(T \mal x_{[\varrho,2]})} \setminus p(T \mal x_{[\varrho,2]}) \\
\\
& = & \b{p(T \mal x_{[\sigma_1, 1]})} \setminus p(T \mal x_{[\sigma_1, 1]}).
\end{eqnarray*}

Here ``$\b{ \ \vphantom{a} \ }$'' refers to taking the Zariski
closure. Since $\b{p(T \mal x_{[\sigma_1, 1]})}$ is of dimension one,
it follows that $T \mal p(x_{[\sigma_2, 2]})$ is a point, i.e.,
$p(x_{[\sigma_2, 2]})$ is fixed by $T$.

So we verified that $f$ is injective. Since it is a regular map of
normal varieties, Zariski's Main Theorem yields that $f$ is an open
embedding. In particular, $f(Y) = q(X)$ is open in $Z$ which is a
contradiction. \endproof

\section{Some Convex Geometry}\label{elconv}

Before coming to the investigation of quotients by the action of
a subtorus,  in this section we recall some elementary properties
 of convex cones
that will be needed later on. Let $V$ be a finite--dimensional real
vector space. In Section 7 we will need the following fact:

\begin{lemma}\label{interiorpoints}
  Let $\varrho$ and $\tau$ be convex cones in $V$
  such that $\varrho^{\circ} \cap \tau \ne \emptyset$. Then
  $\tau^{\circ}$ is contained in the relative interior of $\sigma :=
  \conv(\tau \cup \varrho)$.
\end{lemma}

\proof Since $\tau$ is contained in $\sigma$, it suffices to show that
$\tau^{\circ} \cap \sigma^{\circ}$ is non--empty. Choose $v_{1} \in
\varrho^{\circ} \cap \tau$ and $v_{2} \in \tau^{\circ}$. Then $v :=
v_{1} + v_{2}$ lies in $\tau^{\circ}$. We claim that $v \in
\sigma^{\circ}$. In order to check this, we have to show that every
linear form $u$ contained in the dual cone $\sigma^{\vee}$ of $\sigma$
with $u(v) = 0$ vanishes on $\sigma$. So let $u \in \sigma^{\vee} =
\varrho^{\vee} \cap \tau^{\vee}$ with $u(v) = 0$. Then we obtain that
$u(v_{1}) = u(v_{2}) = 0$.  Consequently $u$ vanishes on $\varrho$ and
$\tau$. This yields $u\vert_{\sigma} = 0$. \endproof

\bigskip

Let $\sigma$ denote a convex (polyhedral) cone in $V$. The set of faces of
$\sigma$ is denoted by $\mathfrak{F}(\sigma)$. Note that the smallest
face of $\sigma$ is $\sigma \cap -\sigma$ and hence equals the maximal
linear subspace contained in $\sigma$.

We consider the following situation: Let $W \subset V$ be any linear
subspace and let $P \colon V \to V/W$ denote the projection. Then
$P(\sigma)$ is a cone in $V/W$. We now want to describe the faces of
$P(\sigma)$ in terms of faces of $\sigma$. The first statement is the
following

\begin{remark}\label{Seitenentsprechung}
\begin{enumerate}
\item There is an injective map from $\mathfrak{F}(P(\sigma))$ to
  $\mathfrak{F}(\sigma)$, given by $\t{\tau}\mapsto
  P^{-1}(\t{\tau})\cap \nobreak\sigma$.
\item If $W \subset (\sigma \cap -\sigma)$, then $\tau \mapsto
  P(\tau)$ is a bijective map from $\mathfrak{F}(\sigma)$ to
  $\mathfrak{F}(P(\sigma))$, inverse to the map in i).
$\bemende$
\end{enumerate}
\end{remark}

%\proof i) Let $\t{\tau}$ be a face of $P(\sigma)$. Then $\t{\tau} =
%\t{u}^{\perp} \cap P(\sigma)$ for some linear form $\t{u} \in
%P(\sigma)^{\vee}$. For $u := \t{u} \circ P$ we have $u \in
%\sigma^{\vee}$ and
%$$
%P^{-1}(\t{\tau}) \cap \sigma = P^{-1}(\t{u}^{\perp}) \cap
%P^{-1}(P(\sigma))\cap\sigma = P^{-1}(\t{u}^{\perp}) \cap \sigma =
%u^{\perp} \cap \sigma \prec \sigma.
%$$
%
%ii) Let $\tau \prec \sigma$. Then there is a linear form $u \in
%\sigma^{\vee}$ with $\tau = u^{\perp} \cap \sigma$. Since by
%assumption, $W$ is contained in the smallest face of $\sigma$, we have
%$W \subset \tau$ and hence $u(W) = 0$.
%
%Consequently there is a linear form $\t{u} \in (V/W)^*$ with $u =
%\t{u} \circ P$. It is clear that $\t{u}\in P(\sigma)^{\vee}$ and
%$\t{u}^{\perp} \cap P(\sigma) = P(\tau)$. So in particular, $P(\tau)$
%is a face of $P(\sigma)$.
%
%Moreover, the maps $\tau \mapsto P(\tau)$ and $\t{\tau} \mapsto
%P^{-1}(\t{\tau})$ are inverse to each other since we have
%$P(P^{-1}(\t{\tau})) = \t{\tau}$ and
%$$
%P^{-1}(P(\tau)) = \tau + W = \tau. \quad \kasten$$

Now let
$\sigma_W$ denote the smallest face of $\sigma$ containing
$W\cap\sigma$. Then $\sigma_W$ is also the largest face of $\sigma$
with $\sigma_W^{\circ} \cap W \ne \emptyset$.

\begin{remark}\label{kleinsteSeite}
  $\h{W} := \sigma_W + W$ is the smallest face of the cone $\sigma +
  W$. In particular, $\h{W}$ is a linear subspace of $V$.
$\bemende$
\end{remark}

%\proof First note that we have $\h{W}^{\vee} = \sigma_W^{\vee} \cap
%W^{\perp}$. We claim that even $\h{W}^{\vee} = \sigma_W^{\perp} \cap
%W^{\perp}$ holds. To see this let $u \in \sigma_W^{\vee} \cap
%W^{\perp}$ and $v \in \sigma_W^{\circ} \cap W$. Then $u(v)=0$ and
%therefore $u(\sigma_W) = 0$ which proves the claim.
%
%So the dual cone of $\h{W}$ is a vector space and hence $\h{W}$ is a
%vector space too. In particular, $\h{W}$ is contained in the smallest
%face of the cone $\sigma + W$, i.e.~in $(\sigma+W)\cap (-\sigma+W)$.
%
%Conversely, let $v\in (\sigma+W)\cap (-\sigma+W)$. Then there are
%$v_1,v_2\in\sigma$ and $w_1,w_2\in W$ such that $v=v_1+w_1=-v_2+w_2$.
%We obtain
%$$v_1+v_2=w_2-w_1\in \sigma\cap W\subset\sigma_W.$$
%That implies
%$v_1,v_2\in\sigma_W$ and hence $v\in\sigma_W+W=\h{W}$.  Altogether we
%have shown that $(\sigma+W)\cap (-\sigma+W)=\h{W}$.  \endproof

Let $\sigma$, $V$, $W$ and $\h{W}$ be as above. We consider the
projections $P^1 \colon V \to V/W$ and $P \colon V \to V/\h{W}$.  As
shown above, $\h{W}$ is the smallest face of $\sigma+W$. Hence Remark
\ref{Seitenentsprechung} yields $1-1$--correspondences
\begin{eqnarray*}
\mathfrak{F}(\sigma +W) \to \mathfrak{F}(P^{1}(\sigma)), & \quad
& \tau \mapsto P^{1}(\tau), \\
\mathfrak{F}(\sigma +W) \to \mathfrak{F}(P(\sigma)), & \quad
& \tau \mapsto P(\tau).
\end{eqnarray*}
In particular, the smallest face of $P^{1}(\sigma)$ is $\h{W} / W$ and
$P(\sigma)$ is strictly convex. Moreover, in these notations we have

\begin{remark}\label{projectingcones}
  For a given face $\tau$ of $\sigma$ the following conditions are
  equivalent:
\begin{enumerate}
\item $(\tau + W) \prec (\sigma + W)$ and $(\tau+W) \cap \sigma
  =\tau$.
\item $P^1(\tau) \prec P^1(\sigma)$ and ${P^1}^{-1}(P^1(\tau)) \cap
  \sigma = \tau$.
\item $P(\tau) \prec P(\sigma)$ and $P^{-1}(P(\tau))\cap \sigma
  =\tau$.  $\bemende$
\end{enumerate}
\end{remark}

%\proof The equivalence of properties i) and ii) is clear.
%``i)$\Rightarrow$iii)'':\enspace We have $P(\tau) = P(\tau + W) \prec
%P(\sigma)$. Moreover, we can conclude from $(\tau + W) \prec (\sigma +
%W)$, that $\tau + W$ contains $\h{W}$ and hence $\tau + W = \tau +
%\h{W}$.  Therefore
%$$P^{-1}(P(\tau)) \cap \sigma = (\tau + \h{W}) \cap \sigma = (\tau+W)
%\cap \sigma =\tau.$$
%
%``iii)$\Rightarrow$i)'':\enspace Note first that $P^{-1}(P(\tau)) \cap
%\sigma = \tau$ means $(\tau + \h{W}) \cap \sigma = \tau$. Since
%$\sigma_W$ is contained in $\h{W}\cap\sigma$, that implies $\sigma_{W}
%\subset \tau$. Hence we obtain $\tau + W = \tau + \h{W}$. Since we
%also have $\sigma + W = \sigma + \h{W}$, Lemma
%\ref{Seitenentsprechung} yields i). \endproof
%

For later use we introduce here the following generalization of the
notion of a fan. Let $N$ denote a lattice and let $\Sigma$ be a finite
set of not necessarily strictly convex cones in $N$. We call $\Sigma$
a {\it quasi--fan}, if $\sigma \in \Sigma$ implies that every face of
$\sigma$ lies in $\Sigma$ and for any two cones $\sigma, \sigma' \in
\Sigma$ the intersection $\sigma \cap \sigma'$ is a face of both,
$\sigma$ and $\sigma'$.

For a given quasi--fan $\Sigma$ in $N$, let $\sigma_{0}$ denote its
minimal element, i.e., $\sigma_0$ is the minimal face of each $\sigma
\in \Sigma$. Consider the primitive sublattice $L := \sigma_0 \cap N$
of $N$ and let $P \colon N \to N/L$ denote the projection. As an
immediate consequence of Remark \ref{Seitenentsprechung} we obtain:

\begin{remark}
  The set $\Delta := \{P_{\RR}(\sigma); \; \sigma \in \Sigma\}$ is a
  fan in $N/L$. $\bemende$
\end{remark}

The concept of a system of fans also has a natural generalization in
this framework: We call a finite family $\S := (\Sigma_{ij})_{i,j \in
  I}$ a system of quasi--fans if it satisfies the conditions
\ref{systemoffans} i) to iii). Again such a system of quasi--fans is
called affine if each $\Sigma_{ii}$ is the quasi--fan of faces of a
single cone $\sigma(i)$.

As in the case of systems of fans, we define a glueing relation on
the set $\mathfrak{F}(\S) := \{(\sigma,i); \; i \in I, \sigma \in
\Sigma_{ii}\}$ of labelled faces of a system $\S$ of quasi--fans and
denote by $\Omega(\S)$ the set of equivalence classes.

A map of two systems $\S$, $\S'$ of quasi--fans in lattices $N$, $N'$
respectively, is pair $(F,\mathfrak{f})$, where $F \colon N \to N'$ is
a lattice homomorphism and $\mathfrak{f} \colon \Omega(\S) \to
\Omega(\S')$ is a map that satisfies the conditions of
\ref{mapofsystemsoffans}.

For practical purposes we note that a map of the systems $\S$ and
$\S'$ of (quasi--) fans is in certain cases induced by a lattice
homomorphism together with a compatible map of the index sets $I$ and
$I'$.

\begin{lemma}\label{mapofsystems of fans2}
  Let $F \colon N \to N'$ be a lattice homomorphism and let $\mu
  \colon I \to I'$, $i \mapsto i'$ be a map such that for any two $i,j
  \in I$ we have
\begin{itemize}
\item[$(*)$] for every $\sigma \in \Sigma_{ij}$ there is a $\sigma'
  \in \Sigma'_{i' j'}$ with $F_{\RR}(\sigma) \subset \sigma'$.
\end{itemize}
Then there is a unique map $(F,\mathfrak{f})$ of the systems of
quasi--fans $\S$ and $\S'$ with $\mathfrak{f}([\sigma,i]) \prec
[\sigma', i']$ for all $\sigma,\sigma'$ as in $(*)$.
\end{lemma}

\proof Let $[\sigma,i] \in \Omega(\S)$. Choose $\sigma' \in
\Sigma'_{i'i'}$ as in Condition $(*)$ and let $\sigma''$ denote the
smallest face of $\sigma'$ with $F_{\RR}(\sigma) \subset \sigma''$, in
other words $\sigma''$ is the face of $\sigma'$ with
$F_{\RR}(\sigma^{\circ}) \subset (\sigma'')^{\circ}$. Set
$$\mathfrak{f}([\sigma,i]) := [\sigma'', i'] \prec [\sigma', i'].$$
In order to see that $\mathfrak{f}$ is well defined and order
preserving, let $[\tau,j] \prec [\sigma,i]$. Choose $\tau'' \in
\Sigma'_{j'j'}$ as above. Since $\tau \in \Sigma_{ij}$, Property $(*)$
yields $\tau'' \in \Sigma'_{i'j'}$ and hence we obtain
$$[\tau'',j'] = [\tau'', i'] \prec[\sigma',i'].$$
To obtain uniqueness of the map $(F,\mathfrak{f})$ of systems of
quasi--fans, note that $\mathfrak{f}([\sigma,i]) \prec [\sigma', i']$
readily implies $\mathfrak{f}[\sigma,i] = [\sigma'',i]$, where
$\sigma''$ is the cone in $\Sigma_{i'i'}$ with
$F_{\RR}(\sigma^{\circ}) \subset (\sigma'')^{\circ}$. \endproof

\begin{remark}\label{maps of affine systems}
Every map from an affine system of quasi--fans to an arbitrary
system of quasi--fans arises from a map of the index sets as in
\ref{mapofsystems of fans2}. $\bemende$
\end{remark}

In contrast to this observation maps from general systems of
quasi--fans may not have a description by a map of the index sets as
the following example shows.

\begin{example} Let $\S = (\Delta_{11})$ be a single fan consisting of
  two maximal cones $\sigma_1$ and $\sigma_2$. Let $I'=\{1,2\}$ and
  for $i\in I'$ let $\Delta'_{ii}$ denote the fan of faces of
  $\sigma_i$ and let $\Delta'_{12} = \Delta'_{21}$ be the fan of faces
  of $\sigma_{1} \cap \sigma_{2}$. Then the identity of $N$ defines a
  unique map of systems of fans from $\S$ to $\S' :=
  (\Delta_{ij}')_{i,j\in I'}$. But there is no map from $I = \{1\}$ to
  $I'$ satisfying $(*)$. \quad $\diamondsuit$
\end{example}

Now, let $\S =(\Sigma_{ij})_{i,j \in I}$ be a system of quasi--fans in
a lattice $N$. Denote by $\sigma_0$ the minimal element of some (and
hence all) $\Sigma_{ij}$. As above, set $L := \sigma_0 \cap N$ and
let $P \colon N \to N/L$ be the projection. Set
$$\Delta_{ij} := \{P_{\RR}(\sigma); \; \sigma \in \Sigma_{ij}\}.$$

\begin{remark}\label{quasi2fan}
  $\t{\S} := (\Delta_{ij})_{i,j \in I}$ is a system of fans in
  $N/L$. The system $\t{\S}$ is affine, if and only if  $\S$ is affine.
 Moreover, the map
$$ \mathfrak{p} \colon \Omega(\S) \to \Omega(\t{\S}), \qquad
[\sigma,i] \mapsto[P_{\RR}(\sigma), i]$$
 is an order--preserving  bijection and $(P,\mathfrak{p})$ is a map of
 systems of quasi--fans. Any further map from $\S$ to a system of fans
 factors uniquely through $(P, \mathfrak{p})$. $\bemende$
\end{remark}

\section{Good Prequotients}

Let $G$ be a reductive complex algebraic group. For an algebraic
action of $G$ on a variety, Seshadri introduced the notion of a good
quotient (see \cite{Se}, Def.~1.5). His notion can be carried over to
the category of prevarieties. Let $X$ be a complex algebraic
prevariety and assume that $G$ acts on $X$ by means of a regular map
$G \times X \to X$.

\begin{definition}
  A $G$-invariant regular map $p \colon X \to Y$ onto a prevariety $Y$
  is called a {\it good prequotient} for the action of $G$ on $X$ if:
\begin{enumerate}
\item $p$ is an affine map, i.e., for every affine open subspace $V$
  of $Y$ the open subspace $U := p^{-1}(V)$ of $X$ is affine,
\item ${\cal O}_{Y}$ is the sheaf $(p_* {\cal O}_{X})^G$ of
  invariants, i.e., for every open set $V \subset Y$ we have ${\cal
    O}_{Y}(V) = {\cal O}_{X}(p^{-1}(V))^{G}$.
\end{enumerate}
\end{definition}

Note that if a good prequotient exists and if both, $X$ and $Y$, are 
 separated, then the good prequotient
is nothing but the good quotient. But in general, even if $X$ is 
separated, 
the action of $G$ may admit a good prequotient but no good quotient
(see Example \ref{hyperbolic}).

As in the case of varieties, a good prequotient is obtained by glueing
algebraic quotients of $G$--stable affine charts.

\begin{lemma}\label{goodprequotcrit} 
  A $G$--invariant surjective regular map $p \colon X \to Y$ is a good
  prequotient for the action of $H$ if there is a covering of $Y$ by
  open affine subspaces $V_{i}$, $i \in I$, such that for every $i \in
  I$ we have
\begin{enumerate}
\item $U_{i} := p^{-1}(V_{i})$ is an open affine subspace of $X$,
\item $p\vert_{U_{i}} \colon U_{i} \to V_{i}$ is an algebraic
  quotient for the action of $G$ on $U_{i}$, i.e., $p\vert_{U_{i}}$
  is given by the inclusion $\CC[U_{i}]^{G} \subset \CC[V_{i}]$. \endproof
\end{enumerate}
\end{lemma}

%This lemma is an immediate consequence of the following elementary
%fact on affine maps of prevarieties:
%
%\begin{lemma}\label{affmorphcrit}
%  A regular map $f : X \to Y$ of complex prevarieties is affine if and
%  only if $Y$ admits a covering by open affine subspaces $V_{i}$ such
%  that each $U_{i} := f^{-1}(V_{i})$ is an open affine subspace of
%  $X$.
%\end{lemma}
%
%\proof Assume that $Y$ admits a covering by open affine subspaces
%$V_{i}$ as above. Let $V$ be an affine open subset of $Y$. First we
%show that $U := p^{-1}(V)$ is hausdorff in the $\CC$--topology and
%hence in fact a variety. Let $x, y \in U$ be different points.
%
%If $p(x) = p(y)$, then $x$ and $y$ are contained in one of the affine
%charts $U_{i}$, where they can be separated by disjoint $\CC$--open
%neighbourhoods. If $p(x) \ne p(y)$, then there are disjoint
%$\CC$--open neighbourhoods $V_x$ and $V_y$ of $p(x)$ and $p(y)$
%respectively in $V$. Their preimages $p^{-1}(V_x)$ and $p^{-1}(V_y)$
%separate $x$ and $y$ in $U$.
%
%Now, for every $i$ choose a covering $(V_{ij})_{j \in J_{i}}$ of $V
%\cap V_i$ by affine open subspaces. Then each $U_{ij} :=
%p^{-1}(V_{ij})$ is affine since the restriction $p_{i} \colon U_{i}
%\to V_{i}$ of $p$ to $U_{i}$ is an affine morphism. Since $U$ is a
%variety, we obtain that $p\vert_{U}$ is an affine map (see
%e.g.~\cite{Ht}, p.~128). In particular, $U$ is affine.  \endproof

\begin{definition}
  A $G$--invariant regular map $p \colon X \to Y$ to a complex
  prevariety $Y$ is called a {\it categorical prequotient}, if every
  $G$--invariant regular map from $X$ to a prevariety factors uniquely
  through $p$.
\end{definition}

Note that categorical prequotients are necessarily surjective. In
analogy to the situation of varieties (see \cite{Se}, p.~516) one
concludes:

\begin{proposition}\label{good2cat}
  Every good prequotient  is  a
  categorical prequotient. \endproof
\end{proposition}

%\proof Let $f \colon X \to Z$ be a $G$--invariant regular map of
%prevarieties. From the affine case (see \cite{Kr}, II.3.2, p.~96) we
%deduce that for every $y\in Y$ its fibre $p^{-1}(y)$ contains
%precisely one closed $G$--orbit $G \mal x$ and that
%
%$$p^{-1}(y) = \{ x' \in X; \; G \mal x \subset \b{G \mal x'} \}.$$
%
%Consequently, since $f$ is a continuous map of topological spaces and
%the points of $Z$ are closed, $f$ is constant on the fibres of $p$.
%Hence there is a map of sets $\t{f} \colon Y \to Z$ such that $f =
%\t{f} \circ p$.
%
%From the corresponding statement on algebraic quotients of affine
%varieties (see e.g. \cite{Kr}, II.3.2, p.~96) one easily concludes
%that $Y$ carries the quotient topology with respect to $p$. So we
%obtain that $\t{f}$ is continuous. Finally, since we have ${\cal
%  O}_{Y} = (p_* {\cal O}_{X})^G$, it follows that $\t{f}$ is regular.
%\endproof
%

Now we specialize to the case that $X$ is a toric prevariety with 
acting torus $T$ and we consider a subtorus $H\subset T$.

\begin{corollary}\label{good2toric}
  If $p \colon X \to Y$ is a good prequotient for the action of $H$ on
  $X$ then $Y$ is a toric prevariety and $p$ is a toric morphism.
\end{corollary}

\proof Choose a covering of $Y$ by open affine subspaces $V_{i}$ such
that the conditions i) and ii) of Lemma \ref{goodprequotcrit} are
satisfied. Consider the action of $H$ on $T\times X$ defined by $h
\mal (t,x) := (t, h \mal x)$. The map
$$
q := \id_{T} \times p \colon T \times X \to T \times Y$$
is a good prequotient for this action, since the sets $T \times V_{i}$
satisfy the conditions of Lemma \ref{goodprequotcrit}. By Proposition
\ref{good2cat}, we obtain a commutative diagram of regular maps
$$\matrix{%
  T \times X & \bigtopmapright{} & X &\cr & & & \cr \lmapdown{q} & &
  \rmapdown{p} & \cr & & & \cr T \times Y & \bigtopmapright{} & Y &,
  \cr} $$
where the horizontal arrows indicate regular $T$-actions. Since $Y$ is
a normal prevariety the claim follows. \endproof

\bigskip

In Theorem \ref{equivcat} we showed that every toric prevariety arises
from an affine system of fans. As the main result of this section we
characterize in terms of affine systems of fans, when the action  of a
subtorus on the toric prevariety admits a good prequotient. For the
corresponding statements on toric varieties we refer to \cite{Sw} and
\cite{Hm}.

Let us first recall the description of the good quotient in the affine
case. Consider an affine toric variety $X_{\sigma}$ where $\sigma$ is
a strictly convex cone in the lattice $N$ and let $L$ be the primitive
sublattice of $N$ corresponding to a subtorus $H$ of the acting torus
of $X_{\sigma}$.

Let $\sigma_{L} := \sigma_{L_{\RR}}$, i.e., $\sigma_{L}$ is the largest face of
$\sigma$ with $L_{\RR} \cap \sigma_{L}^{\circ} \ne \emptyset$
(see also Section \ref{elconv}). Let $\h{L} := N \cap (L_{\RR} +
\sigma_{L})$. Denote by $P\colon N \to N / \h{L}$ the projection. Then
$\t{\sigma} :=P_{\RR}(\sigma)$ is a strictly convex cone in $N /
\h{L}$. Moreover, we have (see e.g. \cite{acha}, Example 3.1):

\begin{remark}\label{affquot}
  The toric morphism $X_{\sigma} \to X_{\t{\sigma}}$ associated to $P$
  is the algebraic quotient for the action of $H$ on $X_{\sigma}$.
$\bemende$
\end{remark}

Now we formulate our criterion for the general case. Let
$N$ be a lattice, let $I$ be a finite index set, and let 
$\S=(\Delta_{ij})_{i,j\in I}$ be an affine system of fans in $N$. 
Recall that for every $i$ there is a strictly convex cone $\sigma(i)$ such
that $\Delta_{ii}$ is the fan of faces of $\sigma(i)$.

Let $X_{\S}$
be the toric prevariety associated to $\S$ and 
let $H$ be a subtorus of its acting torus.
Let $L$ denote the (primitive) sublattice of $N$ corresponding to $H$
and let $P^{1} \colon N \to N/L$ denote the projection. With
these notations our result is the following:

\begin{theorem}\label{characterization}
  The action of $H$ on $X_{\S}$ admits a good prequotient if and
  only if for every $i, j \in I$ and every $\tau\in\Delta_{ij}^{\max}$
  the following holds:
\begin{enumerate}
\item $P^{1}_{\RR}(\tau) \prec P^{1}_{\RR}(\sigma(i))$,
\item ${P^{1}_{\RR}}^{-1}(P^{1}_{\RR}(\tau)) \cap \sigma(i) = \tau$.
\end{enumerate}
If these conditions are satisfied then there is a primitive sublattice
$\h{L}$ of $L$ such that $\h{L}_{\RR}=
L_{\RR} + \sigma(i)_{L}$ for all $i$.
\end{theorem}

In the proof of this result we use the following description of affine
toric morphisms. Let $(F, \mathfrak{f})$ be a map of systems of fans from
$\S$ to an affine system of fans $\S' := (\Delta'_{i'j'})_{i',j' \in
  I'}$ in some lattice $N'$.

\begin{lemma}\label{affinetoricmorphisms}
  The toric morphism $f$ is affine if and only if for every $i'\in I'$
  the set
  $$R(i') := \{ [\tau,j] \in {\Omega}(\S); \; \mathfrak{f} ([\tau,j])
  \prec [\sigma'(i'),i'] \}$$
  contains a unique
  maximal element.
\end{lemma}

\proof Let $i'\in I'$. Using  Lemma \ref{orbitclosures} and the Fibre
Formula \ref{fibreformula}, we obtain the following formula for
the preimage of the maximal affine chart $X'_{i'} :=
X_{[\sigma'(i'),i']}$ of $X_{\S'}$:
$$
f^{-1}(X'_{i'})= \bigcup_{\mathfrak{f}([\tau,j]) \in R(i')} T\mal
x_{[\tau,j]} \,.$$
This open subspace of $X$ is an affine variety if
and only if $R(i')$ contains a unique maximal element.
This proves the claim.  \endproof

\bigskip

{\bf Proof of Theorem \ref{characterization}.}\enspace Assume first
that the conditions i) and ii) are valid. For  $i\in I$ set
$\t{\sigma}(i) := P_{\RR}^1(\sigma(i))$. Note that $\t{\sigma}(i)$
equals $P_{\RR}^1(\sigma(i) + L_{\RR})$. Thus, by
Remark~\ref{Seitenentsprechung}, the smallest face of $\t{\sigma}(i)$
is
$$P_{\RR}^1(\sigma(i)_L + L_{\RR}) = P_{\RR}^1(\sigma(i)_L).$$
Let ${\Sigma}_{ii}$ denote the quasi--fan of faces of  $\t{\sigma}(i)$. 
For $i, j \in I$ define ${\Sigma}_{ij}$ to be
the set of all faces of the cones $P_{\RR}^1(\tau)$, $\tau \in
\Delta_{ij}^{\max}$. Then it follows from Condition i) that
${\Sigma}_{ij}$ is in fact a sub--quasi--fan of ${\Sigma}_{ii}$.

We claim that $({\Sigma}_{ij})_{i,j\in I}$ is a system
of quasi--fans. To show this we need to verify 
${\Sigma}_{ij} \cap {\Sigma}_{jk} \subset
{\Sigma}_{ik}$. Suppose that $\t{\varrho} \in {\Sigma}_{ij} \cap
{\Sigma}_{jk}$, i.e., $\t{\varrho}$ is a face of a cone $\t{\tau} :=
P_{\RR}^1(\tau) \cap P_{\RR}^1(\tau')$ with some $\tau \in
\Delta^{\max}_{ij}$ and $\tau' \in \Delta^{\max}_{jk}$. By Condition
ii),  we have 
$$(P_{\RR}^{1})^{-1}(\t{\tau}) \cap \sigma(j) = \tau \cap \tau',$$
and in particular, $\t{\tau}=P_{\RR}^1(\tau\cap\tau')$.  Since
$\tau \cap \tau'$ lies in $\Delta_{ij} \cap \Delta_{jk} \subset
\Delta_{ik}$ and since $\t{\tau}$ is a common face of $\t{\sigma}(i)$
and $\t{\sigma}(k)$, we can conclude that $\t{\tau} \in {\Sigma}_{ik}$
and hence $\t{\varrho}\in{\Sigma}_{ik}$.

By construction, the maps $P^1$ and $\mu := \id_I$ satisfy the
assumptions of Lemma~\ref{mapofsystems of fans2}. Hence they determine a
unique map of systems of quasi--fans $(P^1,\mathfrak{p^1})$ from 
$\S$ to $({\Sigma}_{ij})_{i,j\in I}$
with $\mathfrak{p^1}([\sigma(i),i]) = [\t{\sigma}(i),i]$.

Since in a system of quasi--fans the minimal elements of the
$\Sigma_{ij}$ all coincide, Remark~\ref{Seitenentsprechung} yields
$\sigma(i)_{L} + L_{\RR} = \sigma(j)_{L} + L_{\RR}$ for all $i,j$. 
So there is a primitive sublattice $\h{L}$ of $L$ such that $\h{L}_{\RR}
= L_{\RR} + \sigma(i)_{L}$ for all $i \in I$. 

Let $Q\colon N/L \to \t{N} := N/\h{L}$ denote the
projection. According to Remark~\ref{quasi2fan}
the sets $\t{\Delta}_{ij} := \{Q_{\RR}(\tau); \; \tau \in \Sigma_{ij}\}$
form a system $\t{\S}$ of fans in $\t{N}$ and the map
$(Q,{\mathfrak q})$ with ${\mathfrak q} \colon 
[\tau,i] \mapsto[Q_{\RR}(\tau),i]$ is universal with respect to maps to systems
of fans.

Let $p \colon X_{\S} \to X_{\t{\S}}$ denote the toric morphism
associated to $(Q\circ P^1,\mathfrak{q}\circ\mathfrak{p})$. 
Since 
the conditions of Lemma~\ref{affinetoricmorphisms} are
satisfied, the morphism $p$
is affine. Moreover, by Remark~\ref{affquot}, for every $i\in I$ the
restriction $p\vert_{X_{i}} \colon X_{i} \to \t{X}_{i}$ is the
algebraic quotient for the action of $H$. Now it follows from Lemma~\ref{goodprequotcrit} that $p$ is a good prequotient for the action of
$H$ on $X$.

Conversely, let $p \colon X \to \t{X}$ be a good prequotient for the action
of $H$ on $X$. By Corollary~\ref{good2toric} and Theorem~\ref{equivcat}, we may assume that $p$ arises from a map
$(P,\mathfrak{p})$ of affine systems of fans $\S$ in $N$ and $\t{\S}$
in $\t{N}$.

The restriction $p_1\colon p^{-1}(\t{T})\to \t{T}$ of $p$ is an
algebraic quotient of affine toric varieties.  So, since $P$ is the
lattice homomorphism associated to $p_{1}$, Remark~\ref{affquot}
implies that $P$ is surjective. Therefore, setting $\h{L} := \ker(P)$,
we can assume that $\t{N} = N / \h{L}$ and $P$ is the canonical
projection.

Since $p$ is an affine surjective toric morphism, we can assume by
Lemma~\ref{affinetoricmorphisms} and the fibre formula that $I =
\t{I}$ and $\mathfrak{p}([\sigma(i), i]) = [\t{\sigma}(i), i]$ with
$\t{\sigma}(i) = P_{\RR}(\sigma(i))$ hold. Moreover, since the $X_{[\sigma(i),i]}$ are
maximal $T$--stable affine open subspaces of $X$, we have
$$ X_{[\sigma(i),i]} = p^{-1}(\t{X}_{[\t{\sigma}(i),i]}).$$
Since the restriction of $p$ to $X_{[\sigma(i),i]}$ is an algebraic quotient for the action of $H$,
Remark \ref{affquot} yields $\sigma(i)_{L} + L_{\RR} = \sigma(j)_{L} +
L_{\RR}$ any two $i, j \in I$.

Now consider $\tau \in \Delta_{ij}^{\max}$ for some $i,j\in I$.  Since
$(P,\mathfrak{p})$ is a map of systems of fans, there is a cone $\t{\tau} \in
\t{\Delta}_{ij}$ such that $\mathfrak{p}([\tau,i]) = [\t{\tau},i]$.  On
the other hand, $\t{\tau} \prec P_{\RR}({\sigma}(i))$ implies by Remark
\ref{Seitenentsprechung} that $\sigma := P_{\RR}^{-1}(\t{\tau}) \cap
\sigma(i)$ is a face of $\sigma(i)$. Clearly we have $\tau \prec \sigma$.

Since we have $\mathfrak{p}([\sigma, i]) = [\t{\tau}, i]\prec
[\t{\sigma}(j),j]$, Lemma \ref{affinetoricmorphisms} yields
$[\sigma,i]\prec[\sigma(j),j]$ and hence $\sigma \in \Delta_{ij}$.
That implies $\tau = \sigma$ and we obtain that $P_{\RR}(\tau) =
\t{\tau}$ and $P_{\RR}^{-1}(P_{\RR}(\tau)) \cap \sigma(i) = \tau$. As
a consequence of Remark \ref{projectingcones} we get conditions i) and
ii). \endproof

\begin{corollary}\label{criterion for good prequotient}
  Let $p\colon X_{\S} \to X_{\S'}$ be a surjective affine
  toric morphism of prevarieties. If the homomorphism of the acting
  tori associated to $p$ has a connected kernel $H$, then $p$ is a
  good prequotient for the action of $H$ on $X_{\S}$.
\end{corollary}

\proof We may assume that $p$ arises from a map $(P, \mathfrak{p})$ of
systems of fans. Since $H$ was assumed to be connected, 
$P$ is surjective and
hence a projection. By Lemma \ref{affinetoricmorphisms} we can assume
that $P_{\RR}(\sigma(i)) = \sigma'(i)$ and
$\mathfrak{p}([\sigma(i),i]) = [\sigma'(i),i]$.

Let $\tau \in \Delta_{ij}^{\max}$ for some $i \ne j$.  Then there is a
cone $\tau' \in \Delta'_{ij}$ such that
$\mathfrak{p}([\tau,i])=[\tau',j]$.  Since $\tau' \prec \sigma'(i)$,
we have $\sigma := P_{\RR}^{-1}(\tau') \cap \sigma(i) \prec \sigma(i)$
and $\mathfrak{p}([\sigma,i])=[\tau',j]\prec [\sigma'(j),j]$. By Lemma
\ref{affinetoricmorphisms}, $\sigma \in \Delta_{ij}$ and hence $\sigma
= \tau$. That proves the claim.  \endproof

\bigskip

If a toric variety $X$ admits a good quotient $p \colon X \to Y$ for the
action of some subtorus $H$ then by definition, $p$ is also a good
prequotient for the action of $H$. The converse of this statement does
not hold, as we see in the following simple example.

\begin{example} \label{hyperbolic}
  The toric variety $X := \CC^2\backslash\{0\}$ is described by the affine
  system $\S = (\Delta_{ij})$ of fans
  in $\ZZ^2$, where $\Delta_{ii}$ for $i=1,2$ denotes the fan of
 faces of $\sigma(1) := \RR_{\ge 0}
  e_{1}$ and $\sigma(2) := \RR_{\ge 0} e_{2}$ and $\Delta_{12} = \{\{0\}\}$.
Consider the subtorus $H := \{(t,t^{-1}); \; t \in \CC^{*}\}$ of the
acting torus $(\CC^{*})^{2}$ of $X$. Then $H$ corresponds to the
sublattice $L$ in $\ZZ^2$ generated by $e_1-e_2$.

The projection $P\colon \ZZ^2\to \ZZ$,
$(x,y)\mapsto x+y$ defines an $H$-invariant toric morphism $p$ from
$X$ onto the complex line with doubled zero with $p(z,w)=zw$ for
$z,w\ne 0$.
%
%\begin{center}
%  \input{hyperbolic.pstex_t}
%\end{center}
%
The morphism $p$ is a good prequotient for the action of $H$, but
there is no good quotient for the action of $H$ (see e.g.  \cite{Sw}
or \cite{Hm}). \quad $\diamondsuit$
\end{example}

We conclude this section with two further examples, showing that both
conditions  of Theorem \ref{characterization} are actually needed.

\begin{example}\label{ifail}
  Consider as in Example \ref{hyperbolic} the affine system of fans $\S$ in
  $\ZZ^2$ defining the toric variety $X := \CC^{2} \setminus \{0\}$.
  Let $L := \RR e_{1}$. Then the subtorus $H$ corresponding to $L$
  equals $\CC^{*}\times \{1\}$.
%\begin{center}
%  \input{ifail.pstex_t}
%\end{center}
The associated projection of lattices is $P^{1} \colon \RR^{2} \to
\RR$, $(x,y) \mapsto y$. Property \ref{characterization} i) is valid
but \ref{characterization} ii) is not. And indeed, the toric morphism
$\CC^2\backslash\{0\}\to \CC$, $(z,w) \mapsto w$ associated to $P^{1}$
is a toric prequotient but not a good quotient.  \quad $\diamondsuit$
\end{example}

\begin{example}\label{iifail}
  Let $\S$ be the affine system of fans in $\ZZ^{6}$ obtained from
$$\sigma(1) := \cone(e_{1}, \ldots, e_{4}), \qquad \sigma_{2} := \cone
(e_{1}, e_{4}, e_{5}, e_{6})$$
 by defining $\Delta_{ii} := {\mathfrak F}(\sigma(i))$
  and $\Delta_{12}:={\mathfrak F}(\sigma(1) \cap \sigma(2))$.
  Define a projection $P^{1} \colon \ZZ^{6} \to \ZZ^{3}$ by
  $$
  \begin{array}{lll} P^{1}(e_{1}) := e_{1}, & P^{1}(e_{2}) := e_{1} +
    e_{3}, & P^{1}(e_{3})
    := -e_{2}, \\
    P^{1}(e_{4}) := e_{1}- e_{3}, & P^{1}(e_{5}) := e_{2}, &
    P^{1}(e_{6}) := e_{1} + e_{3} .
\end{array}$$
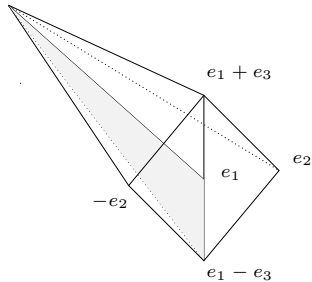
\begin{center}
  \input{iifail.pstex_t}
\end{center}
Let $L := \ker(P^{1})$. Note that
$$P^{1}_{\RR}(\sigma(1)) = \cone(e_{1}+e_{3},
-e_{2}, e_{1}-e_{3}), \qquad P^{1}_{\RR}(\sigma(2)) =
\cone(e_{1}+e_{3}, e_{2}, e_{1}-e_{3}).$$
Thus we see that for the face
$\tau := \cone (e_{1},e_{4}) \in \Delta_{12}^{\max}$ Property
\ref{characterization} i) is not valid. However, Property
\ref{characterization} ii) holds. \quad $\diamondsuit$
\end{example}

\section{The Toric Prequotient}

For actions of subtori of the acting torus of a toric variety, we
introduced in \cite{acha} the notion of a toric quotient. The
analogous concept for the action of a subtorus $H$ of the acting torus
of a toric prevariety $X$ is the following:

\begin{definition}\label{torquot}
  An $H$--invariant toric morphism $p \colon X \to Y$ to a toric
  prevariety $Y$ is called a {\it toric prequotient} for the action of
  $H$ on $X$ if for every $H$-invariant toric morphism $f$ from $X$ to
  a toric prevariety $Z$ there is a unique toric morphism $\t{f}
  \colon Y \to Z$ such that $f = \t{f} \circ p$.
\end{definition}

If $p \colon X \to Y$ is a toric prequotient for the action of a
subtorus $H$ of the acting torus of $X$, then the toric prevariety $Y$
is unique up to isomorphy and will also be denoted by $X \tpq H$. As a
consequence of Proposition \ref{good2cat} and Corollary
\ref{good2toric}, every good prequotient is a toric prequotient. The
aim of this section is to give a constructive proof for the following

\begin{theorem}\label{torquotex}
  Every subtorus action on a toric prevariety admits a toric
  prequotient.
\end{theorem}

In view of Theorem \ref{equivcat}, we prove this result in terms of
affine systems of fans. For the translation of the universal property
of the toric prequotient into the language of systems of fans we
observe:

\begin{remark}
  Let $(F,\mathfrak{f})$ be a map of systems of fans $\S$, $\S'$ in
  lattices $N$, $N'$ respectively, and let $H$ be a subtorus of the
  acting torus $T$ of $X_{\S}$. Then the toric morphism $f \colon
  X_{\S} \to X_{\S'}$ determined by $(F,\mathfrak{f})$ is
  $H$--invariant if and only if the sublattice $L \subset N$
  corresponding to $H$ is contained in $\ker(F)$. $\bemende$
\end{remark}

Now let $N$ be a lattice and let $\S$ be an affine system of
quasi--fans in $N$. Moreover, let $L$ be a primitive sublattice of
$N$. Then the analogue of Definition \ref{torquot} is the following:

\begin{definition}\label{defprequotfan}
  A {\it prequotient} for $\S$ by $L$ is a map of systems of
  quasi--fans $(P,\mathfrak{p})$ from $\S$ to an affine system
  $\t{\S}$ of quasi--fans in a lattice $\t{N}$ such that:
\begin{enumerate}
\item $L \subset \ker(P)$.
\item For every map $(F,\mathfrak{f})$ from $\S$ to an affine system
  of quasi--fans $\S'$ with $F\vert_{L} = 0$, there is a unique map
  $(\t{F},\t{\mathfrak{f}})$ of the systems of quasi--fans $\t{\S}$
  and $\S'$ such that $(F,\mathfrak{f}) = (\t{F},\t{\mathfrak{f}})
  \circ (P, \mathfrak{p})$.
\end{enumerate}
\end{definition}

By Remark \ref{quasi2fan}, for every affine system $\t{S}$ of
quasi--fans in a lattice $\t{N}$ we have  a map that is universal with respect
to maps from $\t{S}$ to affine systems of fans. Thus Theorem
\ref{torquotex} follows directly from Theorem \ref{equivcat} and
the following

\begin{theorem}\label{prequotcalc}
  There is an algorithm to construct for a given affine system of
  quasi--fans $\S$ in $N$ and a primitive sublattice $L$ of $N$ the
  prequotient of $\S$ by $L$.
\end{theorem}

For the proof of this theorem we introduce the following notion. Let
$I$ be a finite index set. We call a collection $\mathfrak{S} :=
(S_{ij})_{i,j\in I}$ of finite sets of cones in $N$ a {\it system of
  related cones} in $N$ if the following conditions are satisfied:

\begin{enumerate}
\item $S_{ii}$ contains precisely one maximal cone $\sigma(i)$,
\item $S_{ij} = S_{ji}$ for all $i, j \in I$,
\item $\tau \in S_{ij}$ implies $\tau \subset \sigma(i) \cap
  \sigma(j)$,
\item If $\tau\in S_{ij}$ then $S_{ij}$ also contains all the faces of
  $\tau$.
\end{enumerate}

A {\it map of two systems of related cones} $\mathfrak{S}$ and
$\mathfrak{S}'$ in lattices $N$ and $N'$ respectively is a pair
$(F,\mu)$, where $F \colon N \to N'$ is a lattice homomorphism and
$\mu \colon I \to I'$, $i\mapsto i'$ is a map of the index sets of
$\mathfrak{S}$ and $\mathfrak{S}'$ such that
\begin{description}
\item{($*$)} for every $\tau \in S_{ij}$ there is a $\tau' \in
  S'_{i'j'}$ with $F_{\RR}(\tau) \subset \tau' $.
\end{description}

Note that every affine system $\S$ of quasi--fans in $N$ is a system
of related cones in $N$. For two affine systems of quasi--fans $\S$ and
$\S'$ in $N$ and $N'$ respectively, any map $(F,\mu)$ from $\S$ to
$\S'$ as map of systems of related cones uniquely determines a map
$(F,\mathfrak{f})$ from $\S$ to $\S'$ as map of systems of fans such that
$$\mathfrak{f}([\sigma,i]) \prec [\sigma'(i'),i']$$
holds for all
$[\sigma,i] \in \Omega(\S)$ (see Lemma \ref{mapofsystems of fans2})
and every map of affine systems of quasi--fans arises in this way.
But a given $(F,\mathfrak{f})$ can arise from different maps of the
systems $\S$ and $\S'$ of related cones.

\bigskip

{\bf Proof of Theorem \ref{prequotcalc}.}\enspace Let $\S =
(\Sigma_{ij})_{i,j \in I}$ be an affine system of quasi--fans in $N$
and let $L$ be a primitive sublattice of $N$.  We use the following
procedure for the calculation of the prequotient
of $\S$ by $L$:

\bigskip

{\it Initialization:}\enspace Set $\t{N} := N / L$ and let $P \colon N
\to \t{N}$ denote the projection. For every $i \in I$ set $\tau^{1}(i)
:= P_{\RR}(\sigma(i))$. For $i,j\in I$ let $S^{1}_{ij}$ denote the set
of faces of the cones $P_{\RR}(\varrho)$, $\varrho \in
\Sigma_{ij}^{\max}$. Set $\mathfrak{S}^{1} := (S^{1}_{ij})_{i, j \in
  I}$.

\bigskip

{\it Loop 1:}\enspace While there are $i,j \in I$, $\varrho \in
{S^{1}_{ij}}^{\max}$, where ``max'' refers to the face relation, with
$\varrho \not\prec \tau^{1}(i)$ do the following: Let $\varrho_{i}$
denote the face of $\tau^{1}(i)$ with $\varrho^{\circ} \subset
\varrho_{i}^{\circ}$. Replace $\tau^{1}(j)$ by $\conv(\tau^{1}(j) \cup
\varrho_{i})$ and replace $S^1_{jj}$ by the set of faces of
$\conv(\tau^{1}(j)\cup \varrho_{i})$.  Remove $\{\varrho' ; \;
\varrho' \prec \varrho \}$ from $S^{1}_{ij}$, $S^{1}_{ji}$ and add
instead $\{\varrho' ; \; \varrho' \prec \varrho_{i} \}$.

\bigskip

{\it Loop 2:}\enspace While there are $i,j,k \in I$ and $\varrho \in
S^{1}_{ij} \cap S^{1}_{jk}$ such that $\varrho \not\in S^{1}_{ik}$,
replace $S^{1}_{ki}$ and $S^{1}_{ik}$ by $S^{1}_{ik} \cup \{\varrho' ;
\varrho' \prec \varrho\}$.

\bigskip

{\it Output:}\enspace For every $i,j\in I$ let
$\t{\tau}(i):=\tau^{1}(i)$ and $\t{\Sigma}_{ij}:=S^{1}_{ij}$.  Set
$\t{\S} := (\t{\Sigma}_{ij})_{i, j \in I}$.

\bigskip

In order to check that the output is in fact well--defined, we have to
show that the loops of the algorithm are finite. This is clear for
Loop 2. For Loop~1 we use a similar argument as in \cite{acha}, proof
of Theorem 2.3:

Since for each $i, j \in I$ the number of maximal cones of
$S^{1}_{ij}$ does not increase when carrying out a step of Loop~1, it
stays fixed after finitely many, say $K$, steps of Loop~1. Let $E
\subset N$ be a minimal set of generators for the cones $\sigma(i)$,
$i \in I$. Then in each step after the first $K$ steps the number
$$
\sum_{i,j \in I}\sum_{\tau \in {S^{1}_{ij}}^{\max}} \vert P(E) \cap
\tau \vert$$
is properly enlarged. This can happen only a finite
number of times, i.e., Loop~1 is finite. Thus we obtain that the
outputs are in fact well--defined.

\bigskip

{\it Claim:}\enspace $\t{\S}$ is an affine system of quasi--fans in
$\t{N}$.  Moreover, $(P,\id_{I})$ is a map of systems of related cones
from $\S$ to $\t{\S}$ and hence defines a map $(P,\mathfrak{p})$ of
the systems of quasi--fans $\S$ and $\t{\S}$ such that
$\mathfrak{p}([\sigma,i]) \prec [\t{\tau}(i),i]$ for all
$[\sigma,i]\in \Omega(\S)$. The map $(P,\mathfrak{p})$ is the
prequotient for $\S$ by $L$.

\bigskip

We prove this claim: After leaving Loop~1, every $S_{ii}^{1}$ is the
quasi--fan of faces of $\tau^{1}(i)$ and every $S^{1}_{ij}$ is a
sub--quasi--fan of the quasi--fan of common faces of $\tau^{1}(i)$ and
$\tau^{1}(j)$: $S^1_{ij}\prec S^1_{ii}\cap S^1_{jj}$.  Note that this
property is not affected in Loop~2.

Thus the quasi--fans $\t{\Sigma}_{ij}$ satisfy Properties i) and ii)
 Definition~\ref{systemoffans}. The transitivity axiom iii) is
guaranteed by Loop~2. In other words, $\t{\S}$ is an affine system of
quasi--fans.

By construction, $(P,\id_{I})$ is a map of the systems $\S$ and
$\t{\S}$ of related cones. Hence there is a unique map
$(P,\mathfrak{p})$ of the systems of quasi--fans $\S$ and $\t{\S}$
with
$$\mathfrak{p}([\sigma,i]) \prec [\t{\tau}(i),i]$$
for all $[\sigma,i]
\in \Omega(\S)$. We have to prove that $(P,\mathfrak{p})$ satisfies
the universal property of the prequotient of $\S$ by $L$. So, let
$(F,\mathfrak{f})$ be a map from $\S$ to an affine system
$\S'=(\Sigma'_{ij})_{i,j\in I'}$ of quasi--fans in a lattice $N'$ such
that $L \subset \ker(F)$. Then there is a lattice homomorphism $\t{F}
\colon \t{N} \to N'$ with $F = \t{F} \circ P$.

Now choose a map $\mu \colon I \to I'$, $i\mapsto i'$ such that
$\mathfrak{f}([\sigma(i),i]) \prec [\sigma'(i'),i']$.  Then
$(\t{F},\mu)$ is a map from the system of related cones
$\mathfrak{S}^1$, defined as in the initialization, to $\S'$ such that
$$
(F,\mu) = (\t{F},\mu) \circ (P,\id_{I})\,.$$
We show inductively
that $(\t{F},\mu)$ remains a map of lists of related cones, when
$\mathfrak{S}^1$ is modified in one of the two loops.

Suppose we are in Loop~1 and there are $i, j \in I$ and $\varrho \in
{S^{1}_{ij}}^{\max}$ with $\varrho \not\prec \tau^{1}(i)$. Let
$\varrho_{i}$ denote the smallest face of $\tau^{1}(i)$ containing
$\varrho$. Then $\varrho^{\circ} \subset \varrho_{i}^{\circ}$. By the
induction hypothesis, there is a cone $\varrho' \in \Sigma'_{i'j'}$
with $\t{F}_{\RR}(\varrho) \subset \varrho'$.

We claim that we even have $\t{F}_{\RR}(\varrho_i) \subset \varrho'$.
To see this note that $\t{F}_{\RR}(\varrho_i) \subset
\t{F}_{\RR}(\tau^{1}(i)) \subset \sigma'(i')$. So there is a face
$\sigma'$ of $\sigma'(i')$ such that $\t{F}_{\RR}(\varrho_i^{\circ})
\subset {\sigma'}^{\circ}$.  On the other hand,
$\t{F}_{\RR}(\varrho^{\circ}) \subset \varrho'\cap
{\sigma'}^{\circ}\ne \emptyset$. That implies $\sigma'\prec\rho'$ and
the claim follows.

Consequently, we obtain $\t{F}_{\RR}(\tau^{1}(j)\cup\varrho_i) \subset
\sigma'(j')$. So the compatibility condition ($*$) for $(\t{F},\mu)$
remains true after replacing $\tau^{1}(j)$ by $\conv(\tau^{1}(j) \cup
\varrho_{j})$, $S^1_{jj}$ by the set of faces of $\conv(\tau^{1}(j)
\cup \varrho_{j})$ and, in $S^{1}_{ij}$, $S^{1}_{ji}$, the faces of
$\varrho$ by those of $\varrho_i$.

Now consider Loop~2 and suppose that there are $i,j,k \in I$ and
$\varrho \in S^{1}_{ij} \cap S^{1}_{jk}$ such that $\varrho\not\in
S^{1}_{ik}$.  By induction hypothesis, there are cones
$\varrho'\in\Sigma'_{i'j'}$ and $\varrho''\in\Sigma'_{j'k'}$ such that
$\t{F}_{\RR}(\varrho)\subset \varrho'\cap\varrho''$.

Since $\varrho'$ and $\varrho''$ are faces of $\sigma'(j')$ they
intersect in a common face and in particular, $\varrho'\cap\varrho''$
lies in $\Sigma'_{i'j'} \cap \Sigma'_{j'k'}$ and hence in
$\Sigma'_{i'k'}$. This shows that ($*$) for $(\t{F},\mu)$ remains true
after adding $\varrho$ and all its faces to $S^{1}_{ik}$ and
$S^{1}_{ki}$.

In other words the pair $(\t{F}, \mu)$ is a map of the systems
$\t{\S}$ and $\S'$ of related cones. Moreover, by definition we have
$$(F,\mu) = (\t{F},\mu) \circ (P,\id_{I})\,.$$
Consequently the associated map of prefans $(\t{F}, \t{\mathfrak{f}})$
from $\t{\S}$ to $\S'$ is a factorization of $(F,\mathfrak{f})$
through $(P,\mathfrak{p})$. We have to check that $(\t{F},
\t{\mathfrak{f}})$ is uniquely determined by this property.

Since $P$ is surjective, $\t{F}$ is determined by $F = \t{F} \circ P$.
Note that, before entering Loop~1 for each $i \in I$, the image
$P_{\RR}(\sigma(i)^{\circ})$ is contained in $\tau^{1}(i)^{\circ}$. By
Lemma \ref{interiorpoints} this property remains valid after enlarging
the cones $\tau^{1}(i)$ as in Loop~1, i.e., we have in fact
$$
\mathfrak{p}([\sigma(i), i]) = [\t{\tau}(i), i]$$
holds for every
$i \in I$. Consequently one obtains $\t{\mathfrak{f}}([\t{\tau}(i),
i]) = \mathfrak{f}([\sigma(i), i])$ for each $i \in I$. Since
$\t{\mathfrak{f}}$ is order--preserving, it is already determined by
this property.  \endproof

\begin{remark}
  The toric prequotient for the action of $H$ on $X_{\S}$ is good if
  and only if the algorithm for constructing the prequotient of $\S$ by
  the sublattice $L$ corresponding to $H$ already terminates after the
  initialization and Property \ref{characterization} ii) holds.
  $\bemende$
\end{remark}

\begin{example}
  Let $\S$ be the affine system of fans in $\ZZ^{5}$ with $\sigma(1)
  := \cone(e_{1}, \ldots, e_{4})$ and $\sigma(2) := \cone(e_{3},
  e_{4}, e_{5})$ and the maximal glueing relation. Define a lattice
  homomorphism $P \colon \ZZ^{5} \to \ZZ^{3}$ by $P(e_{1}) := v_{i}$
  where the vectors $v_{i}$ are situated as indicated below.

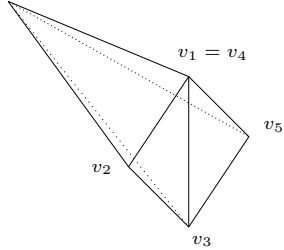
\begin{center}
  \input{termngood.pstex_t}
\end{center}

Clearly we may arrange the $v_{i}$ in such a manner that $P$ is
surjective.  Then for $\tau := \cone(e_{3},e_{4}) \in
\Delta_{12}^{\max}$ we obtain
$${P_{\RR}}^{-1}(P_{\RR}(\tau)) \cap \sigma(1) = \cone(e_{1}, e_{3},
e_{4})$$
and consequently Property \ref{characterization} ii) is not
valid. However, in this situation, the algorithm terminates after the
initialization. \quad $\diamondsuit$
\end{example}

As Example \ref{hyperbolic} indicates, the toric prequotient of a
subtorus action on a toric variety in general differs from its toric
quotient. The two notions are related to each other by the toric
separation (see Section 4):

\begin{remark}\label{intermediatestep}
  For the action of a subtorus $H$ of the acting torus of a toric
  variety $X$, let $p \colon X \to X \tpq H$ be the toric prequotient
  and let $q \colon X \tpq H \to Y$ be the toric separation. Then $q
  \circ p$ is the toric quotient for the action of $H$ on $X$.
  $\bemende$
\end{remark}

In particular, the toric prequotient occurs as an intermediate step in
the construction of the toric quotient. We conclude this section with
an explicit example, showing that both loops of the algorithm are
actually needed.

\begin{example}
  Let us consider the following three cones in $N=\ZZ^7$:
  $$\sigma(1) := \cone(e_1,e_2,e_3), \quad \sigma(2) :=
  \cone(e_2,e_3,e_4,e_5), \quad \sigma(3) := \cone(e_4,e_5,e_6,e_7)
  \,.$$
  Let $\S$ denote the system of fans with these maximal cones
  such that:
  $$\Delta_{12}^{\max} = \sigma(1) \cap \sigma(2), \quad
  \Delta_{23}^{\max} = \sigma(2) \cap \sigma(3), \quad \Delta_{13} =
  \{0\} \,.$$
  Let $P \colon \ZZ^7 \to \ZZ^3$ denote the homomorphism
  given by $P(e_{i}) := v_{i}$, where the $v_{i}$ are vectors in
  $\ZZ^{3}$ situated as in the picture below.

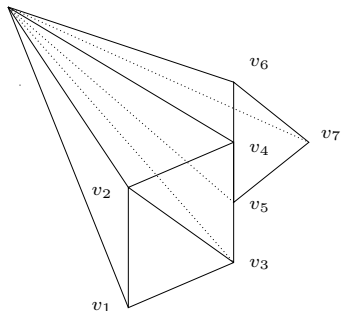
\begin{center}
  \input{loop2.pstex_t}
\end{center}
As before we may assume that $P$ is surjective.  Then after running
through Loop 1 of the prequotient algorithm for $\S$ by $L$ we have
the following cones in $\mathfrak{S}^{1}$:
$$\tau^{1}(1) = P_{\RR}(\sigma(1)), \quad \tau^{1}(2) =
\cone(v_{2},v_{3},v_{6}), \quad \tau^{1}(3) = \cone(v_{3}, v_{6},
v_{7}) \,.$$
and for the families $S^{1}_{ij}$ of subcones we obtain
$$
{S^{1}}^{\max}_{12} = \cone(v_{2}, v_{3}), \quad
{S^{1}}^{\max}_{23} = \cone(v_{3}, v_{6}), \quad {S^{1}}^{\max}_{13} =
\{0\} \,. $$
So we see that $\cone(e_3)$ is contained in $ S_{12} \cap
S_{23}$ but not in $S_{13}$. Consequently the algorithm also enters
Loop~2, where $S_{13}$ is replaced by $\{\cone(e_{3}), \{0\} \}$.
After this the algorithm terminates. \quad $\diamondsuit$
\end{example}

\section{Toric Prevarieties as Prequotients 
of Quasi-Affine Toric Varieties}

In \cite{cox} it is shown that every toric variety occurs as the image
of a good quotient of an open subset of some $\CC^{s}$ by a reductive
abelian group $H$. In fact, a slight modification of Cox's
construction yields that any given toric variety is even the image of a
good toric quotient of an open subset of some $\CC^{s}$ by a subtorus
of $(\CC^{*})^{s}$ (see e.g. \cite{BrVe}, Section 1).

In this section we make related statements in the setting of toric
prevarieties. Let $X_{\S}$ be a toric prevariety arising from a system
of fans $\S = (\Delta_{ij})_{i,j \in I}$ in a lattice $N$. We assume
$\S$ to be affine. According to Theorem \ref{equivcat}, this means no
loss of generality. Moreover, we may assume that $\S$ is irredundant
in the following sense: If $i \ne j$, then $\Delta_{ij}$ is a proper
subfan of $\Delta_{ii}$.

A first aim is to show that $X_{\S}$ occurs as the image of a toric
prequotient of an open toric subvariety of some $\CC^{s}$. Our
construction is the following: For every $i \in I$, let $R_i$ denote
the set of all pairs $(\varrho,i)$, where $\varrho \in
\Delta_{ii}^{(1)}$. Note that $R_{i} \subset \mathfrak{F}(\S)$. Set
$$N' := \bigoplus_{i \in I} \ZZ^{R_{i}}, \qquad \t{N} := N \oplus
N'.$$
For every ray $\varrho\in \bigcup_{i \in I}
\Delta_{ii}^{(1)}$ let $v_{\varrho} \in N$ denote the primitive
lattice vector contained in $\varrho$. Define lattice homomorphisms
$$Q' \colon N' \to N, \quad Q'(e_{(\varrho,i)}) := v_{\varrho}, \qquad
Q := \id_{N} + Q' \colon \t{N} \to N.$$
Here the $e_{(\varrho,i)}$,
$\varrho\in \bigcup_{i \in I} \Delta_{ii}^{(1)}$, denote the canonical
basis vectors of $\ZZ^{R_{i}}$. Now, choose an ordering ``$\le$'' of
$I$ and define an index set $\t{I}$ by
$$
\t{I} := \{ (\tau, i, j); \; i \le j \in I, \; \tau \in
\Delta_{ij}^{\max} \}.$$
So, as a set, $\t{I}$ is isomorphic to the
disjoint union of all $\Delta_{ij}^{\max}$, $i \le j$. For every index
$k = (\tau, i, j) \in \t{I}$ we define a strictly convex lattice cone
$\t{\sigma}_k$ in $N'\subset\t{N}$ by setting
$$
\t{\sigma}_{k} := \cone(e_{(\varrho,l)}; \; l \in \{i,j\}, \;
\varrho \prec \tau ).$$
In particular, if $k = (\sigma, i, i)$ with
the maximal cone $\sigma \in \Delta_{ii}$, then $\t{\sigma}_{k}$ is
the positive quadrant $\RR_{\geq 0}^{R_i}$. By construction, the cones
$\t{\sigma}_{k}$ are the maximal cones of a fan $\t{\Delta}$ in
$\t{N}$, and $X_{\t{\Delta}}$ is isomorphic to an open toric subvariety of
$\CC^{s}$, where $s := \dim(\t{N})$.

Now let $\t{\S}$ denote the affine system of fans determined by
$\t{\Delta}$, i.e., $\t{\S} = (\t{\Delta}_{kk'})_{k,k' \in \t{I}}$,
where $\t{\Delta}_{kk'}$ is the fan of faces of $\t{\sigma}_{k} \cap
\t{\sigma}_{k'}$. 

\begin{lemma}\label{mapmuok}
  For any two elements $k = (\tau,i,j)$ and $k' = (\tau',i',j')$ of
  $\t{I}$ we have
\begin{enumerate}
\item $Q_{\RR}(\t{\sigma}_{k} \cap \t{\sigma}_{k'}) = \{0\} \in
  \Delta_{ii'}$, if $\{i,j\} \cap \{i',j'\} = \emptyset$.
\item $Q_{\RR}(\t{\sigma}_{k} \cap \t{\sigma}_{k'}) = \tau \cap \tau'
  \in \Delta_{ii'}$, if $\{i,j\} \cap \{i',j'\} \ne \emptyset$.
\end{enumerate}
\end{lemma}

\proof Note first that by definition the intersection of
$\t{\sigma}_{k}$ and $\t{\sigma}_{k'}$ is given by
$$
\t{\sigma}_{k} \cap \t{\sigma}_{k'} = \cone(e_{(\varrho,l)}; \; l
\in \{i,j\} \cap \{i',j'\}, \; \varrho \prec \tau, \; \varrho \prec
\tau' ). $$
If $\{i,j\} \cap \{i',j'\}$ is empty, then
$Q_{\RR}(\t{\sigma}_{k} \cap \t{\sigma}_{k'}) = \{0\} \in
\Delta_{ii'}$. So assume that $\{i,j\} \cap \{i',j'\}$ is not empty.
As an example we treat the case $j = j'$. Then $\tau, \tau' \in
\Delta_{jj}$, in particular, $\tau \cap \tau'$ is a face of both,
$\tau$ and $\tau'$. Thus we obtain
$$
Q_{\RR} (\t{\sigma}_{k} \cap \t{\sigma}_{k'}) =
\cone(v_{\varrho}; \; \varrho \prec \tau \cap \tau' ) =
\tau \cap \tau' \in \Delta_{ij} \cap \Delta_{i'j} \subset
\Delta_{ii'}. \quad \kasten $$

\bigskip

By the above lemma, the map 
$\mu \colon \t{I} \to I$, $(\tau,i,j) \mapsto i$,
 satisfies
the condition $(*)$ of Lemma~\ref{mapofsystems of fans2}, and
hence defines a map $(Q,\mathfrak{q})$ of the systems of fans $\t{\S}$
and $\S$. Let $H$ denote the subtorus of the acting torus $\t{T}$ of
$X_{\t{\S}} = X_{\t{\Delta}}$ that corresponds to the primitive
sublattice $\ker(Q) \subset \t{N}$. Then we obtain:

\begin{proposition}\label{precox1cat}
  The toric morphism $q \colon X_{\t{\S}} \to X_{\S}$ associated to
  $(Q,\mathfrak{q})$ is the toric prequotient for the action of $H$ on
  $X_{\t{\S}}$. Moreover, $q$ is even a categorical prequotient for the
$H$--action. 
\end{proposition}

\proof It is clear that $q$ is a surjective toric prequotient.
To see that it is categorical, we have to show the existence of
factorizations. So, let $f \colon X_{\t{\S}} \to Z$ be an $H$-invariant
regular map. As usual, for $i\in I$, let $X_i$ denote the chart
$X_{\Delta_{ii}}$ in $X_{\S}$. For every $k = (\tau,i,j) \in \t{I}$
the lattice homomorphism $Q$ gives rise to a toric morphism
$$q_{k} \; = \; q \vert_{X_{\t{\sigma}_{k}}} \colon X_{\t{\sigma}_{k}} \to
X_{\tau} \;\subset \; X_i\cap X_j \; \subset \; X_{\S}.$$
Note that the
$q_{k}$ are algebraic quotients for the action of $H$ on the
$X_{\t{\sigma}_k}$. In particular, since algebraic quotients are
categorical prequotients (see Proposition \ref{good2cat}), we obtain
for every $k = (\tau,i,j) \in \t{I}$ a regular map $\t{f}_{k} \colon
X_{\tau} \to Z$ such that
$$f \vert_{X_{\t{\sigma}_k}} = \t{f}_{k} \circ q_{k}.$$
We claim that the $\t{f}_{k}$ glue together to a map $\t{f} \colon
X_{\S} \to Z$. To see this, consider first $k=(\tau,i,j)$ and
$k'=(\sigma_i,i,i)$, where $\sigma_i$ denotes the maximal cone in
$\Delta_{ii}$. Recall that
$$\t{\sigma}_{k'} \cap \t{\sigma}_{k} =
\cone(e_{(\rho,i)}; \; \varrho\prec\tau) =: \t{\tau}_{i}.$$
In particular, we have $Q_{\RR}(\t{\tau}_i) = \tau$. Consequently, the
restriction of $q_k$ to $X_{\t{\tau}_i}$ maps $X_{\t{\tau}_i}$ onto
$X_{\tau}$ and is again an algebraic quotient. Since $\t{f}_{k'} \circ
q_k$ coincides with $f$ on $X_{\t{\tau}_i}$, we can conclude that
$$\t{f}_{k'}\vert_{X_{\tau}}=\t{f}_k.$$ 
Thus, to obtain the claim, only the cases $k=(\sigma_i,i,i)$ and
$k'=(\sigma_j,j,j)$ remain to be treated. We have to consider
 $$X_i\cap X_j=\bigcup_{\tau\in\Delta^{\max}_{ij}} X_{\tau}.$$
For every $\tau\in\Delta^{\max}_{ij}$ the previous consideration
yields $\t{f}_k \vert_{X_{\tau}} = \t{f}_{(\tau,i,j)} = \t{f}_{k'}
\vert_{X_{\tau}}$. Therefore $\t{f}_k$ and $\t{f}_{k'}$ in fact coincide on
$X_i\cap X_j$, which proves our claim. By construction we have
$f = \t{f} \circ q$.  \endproof

\begin{corollary}\label{precox1} 
Every toric prevariety $X$ occurs as the image
of a categorical prequotient of an open toric subvariety of some $\CC^s$.
\quad $\kasten$
\end{corollary}

In general, the toric prequotient $q$ constructed above, is not a
good prequotient. In the remaining part of this section we will
 investigate when a given toric prevariety $X_{\S}$ 
can be obtained as a good prequotient
of an open toric subvariety of some $\CC^{s}$. 

Call a toric prevariety
$X$ with acting torus $T$ of {\it affine intersection} if for any two
maximal affine $T$--stable charts $X_{1}, X_{2} \subset X$ their
intersection $X_{1} \cap X_{2}$ is again affine.
Note that every toric variety is of affine intersection while for
toric prevarieties this a proper condition, as the following example
shows:

\begin{example}\label{c2doubled0}
  Let $\sigma:= \cone(e_{1}, e_{2}) \subset \RR^{2}$ and let $\S$
  denote the system of fans with $\Delta_{11} = \Delta_{22}=
{\mathfrak{F}}(\sigma)$  and
$$ \Delta_{12} = \Delta_{21} = \{ \{0\}, \RR_{\ge 0} e_{1}, \RR_{\ge
  0} e_{2} \}. $$
 Then the toric prevariety $X_{\S}$ is just $\CC^{2}$ with doubled
 zero. Clearly $X_{\S}$ is not of affine intersection.
  \quad $\diamondsuit$
\end{example}

\begin{proposition}\label{precox2nec}
  If there is a good prequotient $q \colon \t{X} \to X_{\S}$ with a
  toric variety $\t{X}$, then $X_{\S}$ is of affine intersection.
\end{proposition}

\proof Let $X_{i}$, $i=1,\ldots,r$, be the maximal $T$-stable affine
open subsets of $X_{\S}$ and set $\t{X}_{i} := q^{-1}(X_{i})$.  Then
the restrictions $q_{ij} \colon \t{X}_{i} \cap \t{X}_{j} \to X_{i}
\cap X_{j}$ of $q$ are good prequotients. Since $\t{X}_{i}$ and
$\t{X}_{j}$ are affine, so is $\t{X}_{i} \cap \t{X}_{j}$ and
consequently $X_{i} \cap X_{j}$. \endproof

\bigskip

In the sequel assume that $X_{\S}$ is of affine intersection. Since
$\S$ was assumed to be affine this means  that for any two $i,j
\in I$ the set $\Delta_{ij}^{\max}$ consists of a single cone
$\sigma_{ij}$. We show that $X_{\S}$ occurs as the image of a good
prequotient of an open toric subvariety of some $\CC^{s}$ using the
following generalization of Cox's construction (see \cite{cox}):

Let $R$ denote the set of equivalence classes $[\varrho,i] \in
{\Omega}(\S)$ where $\varrho$ is one--dimensional. Set $N' := \ZZ^R$
and $\t{N} := N \oplus N'$. As before, denote for $\varrho \in \bigcup
\Delta_{ij}^{(1)}$ by $v_{\varrho}$ the primitive lattice vector
contained in $\varrho$. Define lattice homomorphisms
$$
Q' \colon N' \to N, \quad Q(e_{[\varrho,i]}) := v_{\varrho},
\qquad Q := \id_{N} + Q' \colon \t{N} \to N.$$
For every $i\in I$
define a strictly convex cone in $\t{N}$ by setting
$$\t{\sigma}_{i} :=\cone(e_{[\varrho,i]} ; \; \varrho \in
\Delta_{ii}^{(1)})\,$$
Then the cones $\t{\sigma}_{i}$, $i \in I$, are
the maximal cones of a fan $\t{\Delta}$ in $\t{N}$. Let $\t{\S}$ denote
the affine system of fans associated to $\t{\Delta}$.

\begin{lemma}
The homomorphism $Q$ together with $\mu := \id_I$ determines a map
  of systems of fans $(Q,\mathfrak{q})$ from $\t{\S}$ to $\S$.
\end{lemma}

\proof We have to verify condition $(*)$ of Lemma
\ref{mapofsystems of fans2} for $\mu$. Note first that for $i,j \in I$
we have
$$
\t{\sigma}_{i} \cap \t{\sigma}_{k} = \cone (e_{[\varrho,i]}; \;
\varrho \in \Delta_{ij}^{(1)}).$$
Moreover, since $\S$ is affine and
$X_{\S}$ is of affine intersection, $\Delta_{ij}$ is  the fan of
faces of a single cone $\sigma_{ij}$. Hence one obtains condition $(*)$ of Lemma \ref{mapofsystems of fans2} for
$\mu$ from
$$Q_{\RR}(\t{\sigma}_{i} \cap \t{\sigma}_{k}) = \cone(v_{\varrho};
\varrho \in \Delta_{ij}^{(1)})=\sigma_{ij} \in \Delta_{ij}. \qquad
\kasten$$

\bigskip

By construction, the toric morphism $q \colon X_{\t{\S}} \to X_{\S}$
defined by $(Q,\mathfrak{q})$ is surjective and affine. Thus, denoting
by $H$ the kernel of the homomrphism of acting tori associated to $q$,
we infer from Corollary~\ref{criterion for good prequotient}:

\begin{proposition}
  The toric morphism $q \colon X_{\t{\S}} \to X_{\S}$
is a good prequotient for the action of $H$ on  $X_{\t{\S}}$. \endproof
\end{proposition}

Together with Proposition \ref{precox2nec}, the above proposition
yields the following

\begin{theorem}\label{precox2suff}
  For any toric prevariety $X$, the following statements are
  equivalent:
\begin{enumerate}
\item There is an open toric subvariety $U$ of some $\CC^{n}$ and a
  good prequotient $q \colon U \to X$.
\item $X$ is of affine intersection. \endproof
\end{enumerate}
\end{theorem}

\bibliography{}

\end{document}

%% file: nomap.pstex_t
\begin{picture}(0,0)%
\includegraphics{nomap.pstex}%
\end{picture}%
\setlength{\unitlength}{0.00035000in}%
\begingroup\makeatletter\ifx\SetFigFont\undefined
% extract first six characters in \fmtname
\def\x#1#2#3#4#5#6#7\relax{\def\x{#1#2#3#4#5#6}}%
\expandafter\x\fmtname xxxxxx\relax \def\y{splain}%
\ifx\x\y   % LaTeX or SliTeX?
\gdef\SetFigFont#1#2#3{%
  \ifnum #1<17\tiny\else \ifnum #1<20\small\else
  \ifnum #1<24\normalsize\else \ifnum #1<29\large\else
  \ifnum #1<34\Large\else \ifnum #1<41\LARGE\else
     \huge\fi\fi\fi\fi\fi\fi
  \csname #3\endcsname}%
\else
\gdef\SetFigFont#1#2#3{\begingroup
  \count@#1\relax \ifnum 25<\count@\count@25\fi
  \def\x{\endgroup\@setsize\SetFigFont{#2pt}}%
  \expandafter\x
    \csname \romannumeral\the\count@ pt\expandafter\endcsname
    \csname @\romannumeral\the\count@ pt\endcsname
  \csname #3\endcsname}%
\fi
\fi\endgroup
\begin{picture}(8482,2815)(879,-2404)
\put(5356,-556){\makebox(0,0)[lb]{\smash{\SetFigFont{8}{9.6}{rm}$F_{\RR}$}}}
\put(1801,-916){\makebox(0,0)[lb]{\smash{\SetFigFont{8}{9.6}{rm}$\sigma_{1}$}}}
\put(2521,-241){\makebox(0,0)[lb]{\smash{\SetFigFont{8}{9.6}{rm}$\sigma_{2}$}}}
\put(2386,-2356){\makebox(0,0)[lb]{\smash{\SetFigFont{8}{9.6}{rm}$e_{1}$}}}
\put(8461,-1861){\makebox(0,0)[lb]{\smash{\SetFigFont{8}{9.6}{rm}$e_{1}$}}}
\put(9136,-691){\makebox(0,0)[lb]{\smash{\SetFigFont{8}{9.6}{rm}$e_{2}$}}}
\put(3961,-1411){\makebox(0,0)[lb]{\smash{\SetFigFont{8}{9.6}{rm}$e_{2}$}}}
\put(4366, 29){\makebox(0,0)[lb]{\smash{\SetFigFont{8}{9.6}{rm}$e_2 - e_{1} + e_3$}}}
\put(8641,-1411){\makebox(0,0)[lb]{\smash{\SetFigFont{8}{9.6}{rm}$e_{1} + e_2$}}}
\put(9361,-61){\makebox(0,0)[lb]{\smash{\SetFigFont{8}{9.6}{rm}$e_2 - e_{1}$}}}
\end{picture}

%% file: nonsurj.pstex_t
\begin{picture}(0,0)%
\includegraphics{nonsurj.pstex}%
\end{picture}%
\setlength{\unitlength}{0.00035000in}%
\begingroup\makeatletter\ifx\SetFigFont\undefined
% extract first six characters in \fmtname
\def\x#1#2#3#4#5#6#7\relax{\def\x{#1#2#3#4#5#6}}%
\expandafter\x\fmtname xxxxxx\relax \def\y{splain}%
\ifx\x\y   % LaTeX or SliTeX?
\gdef\SetFigFont#1#2#3{%
  \ifnum #1<17\tiny\else \ifnum #1<20\small\else
  \ifnum #1<24\normalsize\else \ifnum #1<29\large\else
  \ifnum #1<34\Large\else \ifnum #1<41\LARGE\else
     \huge\fi\fi\fi\fi\fi\fi
  \csname #3\endcsname}%
\else
\gdef\SetFigFont#1#2#3{\begingroup
  \count@#1\relax \ifnum 25<\count@\count@25\fi
  \def\x{\endgroup\@setsize\SetFigFont{#2pt}}%
  \expandafter\x
    \csname \romannumeral\the\count@ pt\expandafter\endcsname
    \csname @\romannumeral\the\count@ pt\endcsname
  \csname #3\endcsname}%
\fi
\fi\endgroup
\begin{picture}(4332,5175)(1294,-5506)
\put(3376,-5461){\makebox(0,0)[lb]{\smash{\SetFigFont{8}{9.6}{rm}$e_1$}}}
\put(5626,-2536){\makebox(0,0)[lb]{\smash{\SetFigFont{8}{9.6}{rm}$e_2$}}}
\put(2026,-511){\makebox(0,0)[lb]{\smash{\SetFigFont{8}{9.6}{rm}$e_3$}}}
\put(4276,-4336){\makebox(0,0)[lb]{\smash{\SetFigFont{8}{9.6}{rm}$\sigma_1$}}}
\put(2926,-1411){\makebox(0,0)[lb]{\smash{\SetFigFont{8}{9.6}{rm}$\sigma_2$}}}
\end{picture}

%% file: iifail.pstex_t
\begin{picture}(0,0)%
\includegraphics{iifail.pstex}%
\end{picture}%
\setlength{\unitlength}{1579sp}%
\begingroup\makeatletter\ifx\SetFigFont\undefined%
\gdef\SetFigFont#1#2#3#4#5{%
  \reset@font\fontsize{#1}{#2pt}%
  \fontfamily{#3}\fontseries{#4}\fontshape{#5}%
  \selectfont}%
\fi\endgroup%
\begin{picture}(4522,4072)(1104,-3661)
\put(2476,-2536){\makebox(0,0)[lb]{\smash{\SetFigFont{7}{8.4}{\rmdefault}{}{}$-e_{2}$}}}
\put(5626,-1861){\makebox(0,0)[lb]{\smash{\SetFigFont{7}{8.4}{\rmdefault}{}{}$e_{2}$}}}
\put(4276,-511){\makebox(0,0)[lb]{\smash{\SetFigFont{7}{8.4}{\rmdefault}{}{}$e_{1}+e_{3}$}}}
\put(4501,-2086){\makebox(0,0)[lb]{\smash{\SetFigFont{7}{8.4}{\rmdefault}{}{}$e_{1}$}}}
\put(4276,-3661){\makebox(0,0)[lb]{\smash{\SetFigFont{7}{8.4}{\rmdefault}{}{}$e_{1}-e_{3}$}}}
\end{picture}

%% file: termngood.pstex_t
\begin{picture}(0,0)%
\includegraphics{termngood.pstex}%
\end{picture}%
\setlength{\unitlength}{1579sp}%
\begingroup\makeatletter\ifx\SetFigFont\undefined%
\gdef\SetFigFont#1#2#3#4#5{%
  \reset@font\fontsize{#1}{#2pt}%
  \fontfamily{#3}\fontseries{#4}\fontshape{#5}%
  \selectfont}%
\fi\endgroup%
\begin{picture}(4332,3622)(1294,-3661)
\put(2926,-2536){\makebox(0,0)[lb]{\smash{\SetFigFont{7}{8.4}{\rmdefault}{}{}$v_{2}$}}}
\put(4501,-3661){\makebox(0,0)[lb]{\smash{\SetFigFont{7}{8.4}{\rmdefault}{}{}$v_{3}$}}}
\put(4276,-736){\makebox(0,0)[lb]{\smash{\SetFigFont{7}{8.4}{\rmdefault}{}{}$v_{1} = v_{4}$}}}
\put(5626,-1861){\makebox(0,0)[lb]{\smash{\SetFigFont{7}{8.4}{\rmdefault}{}{}$v_{5}$}}}
\end{picture}

%% file: loop2.pstex_t
\begin{picture}(0,0)%
\includegraphics{loop2.pstex}%
\end{picture}%
\setlength{\unitlength}{1579sp}%
\begingroup\makeatletter\ifx\SetFigFont\undefined%
\gdef\SetFigFont#1#2#3#4#5{%
  \reset@font\fontsize{#1}{#2pt}%
  \fontfamily{#3}\fontseries{#4}\fontshape{#5}%
  \selectfont}%
\fi\endgroup%
\begin{picture}(4972,4544)(1104,-4133)
\put(6076,-1411){\makebox(0,0)[lb]{\smash{\SetFigFont{7}{8.4}{\rmdefault}{}{}$v_{7}$}}}
\put(2476,-2311){\makebox(0,0)[lb]{\smash{\SetFigFont{7}{8.4}{\rmdefault}{}{}$v_{2}$}}}
\put(2476,-4111){\makebox(0,0)[lb]{\smash{\SetFigFont{7}{8.4}{\rmdefault}{}{}$v_{1}$}}}
\put(4951,-286){\makebox(0,0)[lb]{\smash{\SetFigFont{7}{8.4}{\rmdefault}{}{}$v_{6}$}}}
\put(4951,-1636){\makebox(0,0)[lb]{\smash{\SetFigFont{7}{8.4}{\rmdefault}{}{}$v_{4}$}}}
\put(4951,-2536){\makebox(0,0)[lb]{\smash{\SetFigFont{7}{8.4}{\rmdefault}{}{}$v_{5}$}}}
\put(4951,-3436){\makebox(0,0)[lb]{\smash{\SetFigFont{7}{8.4}{\rmdefault}{}{}$v_{3}$}}}
\end{picture}

%% file: subtorac.bbl
\begin{thebibliography}{acha}
%
\bibitem%[AC;Ha,1]%
{acha} A.~A'Campo-Neuen, J.~Hausen: Quotients of Toric
  Varieties by the Action of a Subtorus. T\^{o}hoku Math. J. {\bf 51},
  1--12 (1999).
%
\bibitem%[AC;Ha,2]%
{acha3} A.~A'Campo-Neuen, J.~Hausen: Examples and
Counterexamples for Existence of Categorical Quotients. Konstanzer
Schriften in Mathematik und Informatik {\bf 80} (1999),
to appear in Documenta Math..
%
\bibitem%[BB]%
{BB} A. Bia\l ynicki-Birula: Finiteness of the Number of
  Maximal Open Subsets with Good Quotients. Transform. Groups {\bf 3},
No.~4, 301--319 (1998).
%
\bibitem%[Br;Ve]%
{BrVe} M. Brion, M. Vergne: An Equivariant
  Riemann--Roch Theorem for Complete Simplicial Toric
  Varieties. J. reine angew. Math. {\bf 482}, 67--92 (1996).
%
\bibitem%[Co]%
{cox} D. Cox: The Homogeneous Coordinate Ring of a Toric
  Variety. J. Algebraic Geometry {\bf 4}, 17-51 (1995).
%
\bibitem%[Fu]%
{Fu} W.~Fulton: Introduction to Toric Varieties. Princeton
  University Press, Princeton, 1993.
%
\bibitem%[Ha]%
{ha} J.~Hausen: On W\l odarczyk's Embedding Theorem. To
  appear in Int. J. Math.
%
\bibitem%[Hm]%
{Hm} H.~Hamm: Very good quotients of Toric Varieties. 
J.~W.~Bruce et al.  (eds.), Real and Complex Singularities. Proceedings
of the 5th workshop, Sao Carlos,
Brazil, July 27--31, 1998. Boca Raton, FL: Chapman \& Hall/CRC. Chapman
Hall/CRC Res. Notes Math. 412, 61-75 (2000). 
%
%
%\bibitem[Ht]{Ht} R.~Hartshorne: Algebraic Geometry. Springer, Berlin,
%1977.
%
%\bibitem[Kr]{Kr} H. Kraft: Geometrische Methoden der Invariantentheorie. 
%Vieweg, Braunschweig, 1984.%
%
\bibitem%[Ka;St;Ze]%
{KaStZe} M. Kapranov, B. Sturmfels, A. V. Zelevinsky:
  Quotients of Toric Varieties. Math. Ann. {\bf 290}, 643--655 (1991)
%
%\bibitem%[Po;Vi]%
%{PV} V.~L.~Popov, E.~B.~Vinberg: Invariant Theory. 
%In: Algebraic Geometry IV (A.~N.~Parshin, I.~R.~Shafarevich, eds.), 
%Encyclopaedia of Mathematical Sciences {\bf 55}, Springer, Berlin, 1994.
%
\bibitem%[Se]%
{Se} C.S. Seshadri: Quotient Spaces Modulo Reductive
  Algebraic Groups. Ann. of Math. {\bf 95}, 511-556 (1972).
%
\bibitem%[Su]%
{Su} H. Sumihiro: Equivariant completion. Journal of
  Math. Kyoto University {\bf 14}, 1--28 (1974).
%
\bibitem%[Sw]%
{Sw} J.~\'Swi\c{e}cicka: Quotients of Toric Varieties by Actions of Subtori. 
To appear in Colloq. Math.
%
\bibitem%[W\l]%
{Wlodar} J.~W\l odarczyk: Embeddings in toric varieties,
  J.~Algebraic Geometry {\bf 2} (1993), 705--726.
%
\end{thebibliography}
